\documentclass[11pt]{article}
\usepackage{amsmath,amssymb}
\usepackage{latexsym,enumerate}
\usepackage{amsmath,amsthm,amsopn,amstext,amscd,amsfonts,amssymb}
\usepackage{graphicx,natbib}
\usepackage[T1]{fontenc}
\usepackage{amsfonts}

\newtheorem{proposition}{Proposition}
\newtheorem{lemma}{Lemma}
\newtheorem{theorem}{Theorem}
\newtheorem{definition}{Definition}

\title{Deconvolution in white noise with a random blurring function}
\author{Thomas Willer}
\date{}

\begin{document}

\maketitle
\pagenumbering{arabic}
Laboratoire de Probabilit{\'e}s et Mod{\`e}les Al{\'e}atoires, Universit{\'e} Paris
VII (Denis Diderot), 175 rue de Chevaleret, F-75013 Paris, France. \\

\begin{center}
\textbf{Abstract}
\end{center}

We consider the problem of denoising a function observed after a convolution with a random filter independent of the noise and satisfying some mean smoothness condition depending on an ill posedness
coefficient. We establish the minimax rates for the $L^{p}$ risk over balls of periodic Besov spaces with
respect to the level of noise, and we provide an adaptive estimator achieving these rates up to log factors. Simulations were performed to
highlight the effects of the ill posedness and of the distribution of
the filter on the efficiency of the estimator.

\bigskip

\textit{Keywords:} Adaptive estimation; Deconvolution; Inverse problem; Minimax risk; Nonparametric
estimation; Wavelet decomposition.

\section{Motivations and preliminaries}

\subsection{Inverse problems in practice}

Deconvolution is a particularly important case in a more general setting of
problems, known as inverse problems. They consist in recovering an
unknown object $f$ from an observation $h_{n}$ corresponding to $H(f)$
corrupted by a white noise $\xi$, for some operator $H$. The model is of the kind:
\begin{equation}
\label{inverse} h_{n}=H(f)+\sigma n^{-1/2}\xi, \quad \forall n\ge
1.
\end{equation}
Inverse problems appear in many scientific domains.
Several applications can be found for example in \cite{ofta} in
various domains such as meteorology, thermodynamics and mecanics.
Deconvolution, in particular, is a common problem in signal and
image processing (see \cite{Bertero}). It appears notably in light
detection and ranging devices, computing distances to an object
by measuring the lapse of time between the emission of laser
pulses and the detection of the pulses reflected by the object.
In the underlying model $f$ is a distance to an object measured up to small gaussian errors after being blurred by a convolution phenomenon due to the fact that
the system response function of the device is longer than the time
resolution interval of the detector. Several papers deal with this
application of deconvolution methods, for example \cite{Harsdorf}
or \cite{deconvolution}.

\bigskip

In some cases, it is difficult to know \textit{a priori} the underlying
operator which transformed the object to be determined into the
observed data. This problem appears notably when the
operator is sensitive to even slight changes in the experimental
conditions, or is affected by external random effects that cannot
be controlled, and thus changes for every observation. In these
conditions, a framework with a random operator is more adapted
than a setting with a fixed deterministic operator.

As an example
let us consider an inverse problem of reconstruction in a
tomographic imagery system, borrowed from \cite{ofta}. The problem
is to find the density of activity $f$ of a radioactive tracer by
collecting the $\gamma$ photons which it radiates on a detector. The
framework is illustrated on Figure $1$. The setting is such that
only the photons transmitted perpendicularly to the detector are taken
into account. A given pixel $A_d$ of the detector collects a
number of photons that depends on the density of activity $f$
along some segment $[F A_d]$, where $F$ is the focal point towards
which $A_d$ is headed. Each point $M$ of this segment transmits a
contribution $f(M)$ towards $A_d$ but the pixel detects only
$a(M,A_d)f(M)$ photons from $M$ because the radiation diminishes
after it has gone across the fluid between $M$ and $A_d$. So the
following quantity is observed on the pixel $A_d$:

$$X_{\mu}f(F,A_d)=\int_{M \in [F,A_d]}f(M)a(M,A_{d})dM,$$ and the function $a$ can be put in the following form :

$$a(M,A_d)=\exp \bigr[ -\int_{M' \in [M,A_d]} \mu(M')dM'\bigl],$$
where $\mu$ is a coefficient quantifying the radiation fading around $M'$. On figure $1$ several zones
characterized by different densities of activity and different
coefficients $\mu$ are represented. If $\mu$ is constant along the segment $[F
A_d]$, then recovering $f$ is a deconvolution problem.

In practice
the cartography of $\mu$ is not well known $a$ $priori$. There is
a different function for each pixel and this function depends on
the characteristics of the fluid where the tracers were injected.
Complementary measures and reconstruction algorithms are necessary
to obtain it. In this context a probabilistic model is useful,
where $\mu$ is a random function determined $a$
$posteriori$ thanks to additionnal measures.

\begin{figure} \centering
\begin{minipage}[h]{0.6\textwidth}
\centering
\includegraphics[width=3.5in,height=2in]{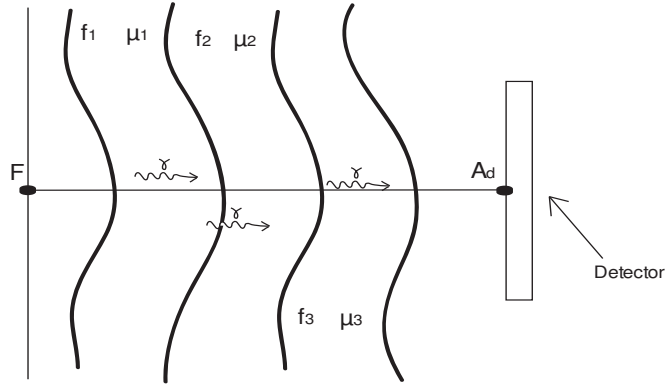}
\caption{\footnotesize Reconstruction of a density of activity}
\end{minipage}
\end{figure}

\subsection{Estimation in inverse problems with random operators}

In the case of deterministic operators, inverse problems have been studied in many papers in a general
framework where \eqref{inverse} holds with some linear
operator $H$. Two main methods of estimation are generally used to
recover $f$ from the observation: singular value decomposition
(SVD) and Galerkin projection methods. The former uses a
decomposition of $f$ on a basis of eigenfunctions of $H^T H$,
which can be hard to perform if $H$ is difficult to diagonalize.
The latter uses a decomposition of $f$ on a fixed basis adapted to
the kind of functions to be estimated and then consists in
solving a finite linear system to recover the coefficients of $f$. Wavelet
decomposition is a very useful tool in such settings, see
\cite{Donoho95} and \cite{Silverman}.

Among others, a method combining
wavelet-vaguelettes decompositions and Galerkin projections can be
found in \cite{marc}, whereas a sharp adaptive SVD estimator can
be found in \cite{cavalier_tsybakov}. Concerning the deconvolution
problem, wavelet-based estimation techniques were developed in
\cite{Pensky}, \cite{Walter}, \cite{Fan}, \cite{Kalifa} and
\cite{deconvolution}. Multidimensional situations have also been
considered: minimax rates and estimation techniques can be found in \cite{tsybakov4}.

\bigskip

Generalisations of inverse problems to the case of random
operators have been made in several recent papers. First, random
operators enable to treat situations where, in practice, the
operator modifying the object to be estimated is not exactly known
because of errors of measure. In such settings, equation
\eqref{inverse} holds with an unknown deterministic operator $H$, and additionnal noisy observations provide a random
operator $H_\delta$ where $\delta$ is a level of noise :
$H_{\delta}=H(f)+\delta \xi.$ The problem is to build an estimator
of $f$ based on the data $(h_{n}, H_{\delta})$ achieving
minimax rates. Several adaptive estimation methods have been
developed in this case. Some are based on SVD methods such as in
\cite{cavalier}, whereas estimators based on Galerkin projection
methods were developed in \cite{efro} or \cite{hoffmann}.

Random operators also appear quite naturally in models
where the evolution of a random process is influenced by its past.
For example let us consider the problem of estimating an unknown function
$f$ thanks to the observation of $X_{n}$ ruled by the following
equation (called stochastic delay differential equation, SDDE in
short):
\begin{align*}
dX_{n}(t) &=(\int_{0}^{r}X_{n}(t-s) f(s)ds) dt+\sigma n^{-1/2}dW(t) \quad \forall t \geq 0,\\
X_{n}(t) &= F(t) \quad \forall t \in [-r,0].
\end{align*}
This problem is close to problem \eqref{eqdep}: a
convolution of the unknown function with the random filter $X_n$
is observed with small errors. However this filter is not
independent from $W$ so our results do not apply to this
particular problem. Numerous estimation results in SDDEs can be
found in \cite{markus} and in \cite{markus2}, with a different asymptotic framework.

\bigskip

The organisation of the paper is as follows. Section $2$, $3$ and $4$ present
respectively the model, the estimator and the main results. Section
$5$ gives simulation results where the behaviour of the estimator is
investigated for several distributions of
the random filter, and section $6$ gives the proofs of the theorems.

\section{The model}

We consider the following deconvolution problem. Let $(\Omega,
\mathcal{A},P )$ be a probability space and $W$ a standard Wiener
process on this space. For a given $n \in \mathbb{N}^{*}$ we
observe the realizations of two processes $X_{n}$ and $Y$ linked
in the following way:

\begin{equation}
  \label{eqdep}
  \begin{cases}
    dX_{n}(t) &= f\star Y(t)dt+\sigma n^{-1/2}dW(t), \quad \forall t \in
    [0,1], \\
    X_{n}(0) &= x_{0},
  \end{cases}
\end{equation} where $\star$ denotes the convolution : $f\star Y(t)=\int_{0}^{1}
f(t-s)Y(s)ds,$ $x_{0}$ is a deterministic initial condition and $\sigma$ is
a positive known constant.

\bigskip

The problem is to estimate the $1$-periodic function $f$
when $Y$ is independent of $W$ and satisfies some condition of
smoothness.

\subsection{The target function}

We introduce functional spaces especially useful to describe the target
functions. For a given $\rho>1$, let us first denote by $L^\rho $ the
following space:
$$L^{\rho}([0,1])=\{ f:\mathbb{R} \mapsto \mathbb{R} \quad|\quad f \mbox{ is }
1-periodic, \quad \mbox{ and }\int_{0}^{1} |f|^{\rho} <
\infty\}.$$
Secondly we use periodic Besov spaces which are defined
thanks to the modulus of continuity in a similar way as in the non
periodic case (see \cite{deconvolution} for the exact definition).
They have the advantage of being very general, including
spatially unsmooth functions, and of being very well suited to wavelet
decompositions. Indeed, the following
characterization holds under
several conditions on the wavelet basis similar to the
conditions in the general case (which can be found in
\cite{wavelets}):
$$B^{s}_{p,q}([0,1])=\{ f \in L^{p}([0,1])\quad| \quad \| f \| _{s,p,q}:=\bigl( \sum_{j \le 0}2^{j(s+1/2-1/p)q}(\sum_{0 \le k \le 2^{j}}
|\beta _{j,k}|^{p})^{q/p} \bigr)^{1/q} < \infty \}.$$

\bigskip

We investigate the maximal
error when $f$ can be any function in a ball of a periodic Besov
space $B^{s}_{p,q}([0,1])$ of radius $R$ and when the estimation error is measured by the $L^{\rho}$-loss. We suppose that $s>\frac {1}{p}$ so that $f$ is continuous an hence its $L^{\rho}$-norm exists.

\begin{definition}
  For given $R > 0$, $p> 1$, $q>1$ and $s>\frac {1}{p}$,
define :
$$M(s,p,q,R)=\{f \in B^{s}_{p,q}([0,1])\quad | \quad \| f \| _{s,p,q} \leq R \}.$$
\end{definition} 
\noindent Our aim is to determine the rate of the following minimax risk for $\rho>1$:
$$R_{n}:=\inf_{\hat {f}_{n}} \sup_{f \in M(s,p,q,R)} E_{f}(\| \hat
{f}_{n}-f \| _{\rho}),$$ where the infimum is taken over all
$\sigma ((X_{n}(t),Y(t))_{t\in[0,1]}))-$measurable estimators
$\hat {f}_{n}$.

\subsection{The filter}

We assume that the blurring function $Y$ is a random process independent of $n$, $f$, and (in probabilistic terms) of the
process $W$, and taking its values in $L^2([0,1])$.

\bigskip

Throughout this paper, we will use the following notations for two
functions $A$ and $B$ depending on parameters $p$ :

\begin{itemize}
\item {$A\lesssim B$ means that there exists a positive constant
$C$ such that for all $p$,$\quad A(p)\le CB(p)$,} \item {
$A\gtrsim B$ means that $B\lesssim A$,} \item { $A\asymp B$ means
that $A\lesssim B$ and $A\gtrsim B$.} \end{itemize}

For $j\in \mathbb{N}$ we introduce two random variables $L^{Y}_{j}$ and $U^{Y}_{j}$ (whenever they exist) linked to the smoothness of the process $Y$:

$$L^{Y}_{j}=\frac {\sum_{l=2^{j}}^{2^{j+1}-1}
|Y_{l}|^{2}}{2^{j}},\quad \mbox{and} \quad U^{Y}_{j}=\frac {\sum_{l=0}^{2^{j+1}-1}
|Y_{l}|^{-2}}{2^{j}},$$ where $(Y_{l})_{l \in \mathbb{Z}}$ are the
Fourier coefficients of $(Y(t))_{t \in [0,1]}$.

\bigskip

To establish the lower (resp upper) bound of the minimax risk, we impose the
following control on the distribution of $L^{Y}_{j}$ (resp $U^{Y}_{j}$), which
implies that the Fourier coefficients are not too large (resp small):

\bigskip

\noindent \textit{$C_{low}$: There exists a constant $\nu\ge 0$ such that, for all
$j\in \mathbb{N}$:
\begin{equation*}
E(L^{Y}_{j}) \lesssim 2^{-2\nu j}.
\end{equation*}} \textit{$C_{up}:\: \forall l\in \mathbb{Z}, Y_{l}\ne
0 \mbox{ almost surely},
\mbox{ and there exist }\nu\ge 0, c>0, \alpha>0 \mbox{ such that, for all }j \in \mathbb{N}:$
$$\forall t\ge 0,\quad P\big(U^{Y}_{j}\ge t 2^{2\nu j}\big) \lesssim e^{-ct^{\alpha}}.$$}

All those conditions are satisfied if the Fourier
Transform $\hat{Y}$ of the process $Y$ has the following form: $|\hat{Y}(w)| =
\frac{T(w)}{(1+w^{2})^{\nu /2}}$, where $T$ is a positive random process
with little probability of taking small or high values (for example bounded almost surely by deterministic
constants). This case
includes for example gamma probability distribution functions with
some random scale parameter, which will be used further. On the contrary, condition $C_{up}$ does not hold for
filters with realizations belonging to supersmooth functions, ie $Y$
such that $|\hat{Y}(w)|=T(w)\frac{e^{-B|w|^{\beta}}}{(1+w^{2})^{\nu
    /2}}$, for some constants $B, \beta > 0$ and with $T$ as
before. Results on deconvolution of supersmooth functions can be found
in \cite{butucea}.

\section{Adaptive estimators}

We first build an adaptive estimator, nearly achieving
the minimax rates exposed in the next section, which is close to the
one developed in \cite{deconvolution} in the
case of a deterministic filter $Y$. The method combines elements of
the SVD methods (deconvolution thanks to the Fourier basis) and of
the projection methods (decomposition on a wavelet basis adapted
to the target functions).

\bigskip

Let us set $R_j=\{0,\dots,2^{j}-1 \}$ for all $j\in \mathbb{N}$, and let $(\Phi_{j,k}, \Psi_{j,k})_{j,k \in \mathbb{Z}}$ denote the
periodized Meyer wavelet basis (see \cite{Meyer} or
\cite{Mallat} for details). For convenience the following notations will be used further: $R_{-1}=\{0\}$ and $\Phi_{-1,0}=\Psi_{0,0}$. Any $1$-periodic target function $f$ belonging to
$M(s,p,q,S)$ has an expansion of the kind:

$$f=\sum_{j\ge -1,\: k \in R_j } \beta_{j,k} \Psi_{j,k},$$ where $$\beta_{j,k}=\int_0 ^1
f\Psi_{j,k}.$$

We estimate $f$ by estimating its wavelet coefficients. Let $(e_l
(t))=(\exp(2\pi ilt))_{l\in \mathbb{Z}}$ denote the Fourier basis, and let $(\Psi_{j,k,l})_{l\in \mathbb{Z}}$, $(f_l)_{l \in
\mathbb{Z}}$ and $(Y_l)_{l \in \mathbb{Z}}$ be the Fourier
coefficients of the functions $\Psi_{j,k}$, $f$ and $Y$. Set also: $W_l=\int _0 ^1 e_l (t) dW(t) $ and $X^n _l =\int _0 ^1 e_l (t)
dX_n (t).$ Then by Plancherel's identity we have:
$$\beta_{j,k}=\sum_{l\in
\mathbb{Z}} f_l \Psi_{j,k,l}.$$ Moreover $\int_0^1(f \star Y) \bar e_l
= f_l Y_l$, so equation \eqref{eqdep} yields:
$$X^{n}_l = f_l Y_l + \sigma n^{-1/2}W_l,$$ and thus if we suppose
that $Y_l \ne 0$ almost surely for all $l$, $f_l$ can naturally be
estimated by $\frac {X^{n}_l}{Y_l}$ and we set:
$$\hat{\beta}_{j,k}=\sum_{l\in
\mathbb{Z}} \frac {X^{n}_l}{Y_l} \Psi_{j,k,l}.$$ Then a
hard thresholding estimator is built with the following values for the
thresholds $\lambda_{j}$ and the highest resolution level $j_1$:
$$2^{j_1}=\{n/ (\log n)^{1+\frac{1}{\alpha}} \}^{1/(1+2\nu)}, $$
$$\lambda_j = \eta 2^{\nu j} \sqrt{(\log n)^{1+\frac{1}{\alpha}}/n },$$ where $\eta$ is
a positive constant larger than a threshold (which is determined in section $6$).

\bigskip

Finally the following estimator achieves the minimax rates up to log factors when the filter satisfies condition $C_{up}$:

\begin{equation}
\label{estimateur}
\hat f_n^{D}=\sum_{(j,k)\in \Lambda_{n}}
\hat{\beta}_{j,k} I_{ \{
  |\hat{\beta}_{j,k}| \geq \lambda_{j}  \} }  \Psi_{j,k},
\end{equation} where $\Lambda_{n}=\{ (j,k)\in \mathbb{Z}^{2}\: | \: j\in
\{-1,\dots,j_1 \}, \: k\in R_j \}$.

\bigskip

Moreover we also introduce a slightly different estimator $\hat f_n^{R}$ with random
thresholds instead of deterministic ones (hence the superscript R
instead of D), ie with $j_1$ and $\lambda_j$ replaced by $j_2$ and $\tau_j$:
$$2^{j_2}=\{ n / \log n \}^{1/(1+2\nu)},$$
$$\tau_j = \eta'  \sqrt{U_j^Y \log n/n },$$ where $\eta'$ is a large
enough constant. The theoretical performances of $\hat f_n^{R}$ will be studied in a
separate publication, here only a simulation study is provided.

\section{Main results}

Let $\rho> 1$, $R>0$, $p> 1$, $q> 1$ and $s> 1/p$. We distinguish three cases for the regularity parameters characterizing the
target functions according to the sign of $\epsilon=\frac{2s+2\nu +1}{\rho} -\frac{2\nu +1}{p}$:

\bigskip

\noindent the sparse case ($\epsilon <0$), the
critical case ($\epsilon= 0$) and the
regular case ($\epsilon>0$).

\bigskip

Let us introduce the two following rates:
\begin{equation*}
r_n(s,\nu)=\big(\frac{1}{n}\big)^{\frac {s} {2s+2\nu +1} } ,\quad s_n(s,p,\rho,\nu)=\big(\frac{\log(n)}{n}\big)^{\frac{s-1/p+1/\rho}{2s+2\nu
+1-2/p}}.
\end{equation*}

\begin{theorem} Under condition $C_{low}$ on $Y$:
\begin{align*}
&r_n(s,\nu)^{-1}R_{n}\gtrsim 1 \quad \mbox{in the regular case,}\\
&s_n(s,p,\rho,\nu)^{-1}R_{n}\gtrsim
1 \quad \mbox{in the sparse and critical cases.}
\end{align*}
\end{theorem}

\begin{theorem} Under condition $C_{up}$ on $Y$:
\begin{align*}
&r_n(s,\nu)^{-1}R_{n}\lesssim 1 \quad \mbox{in the regular case,}\\
&s_n(s,p,\rho,\nu)^{-1}R_{n}\lesssim
1 \quad \mbox{in the sparse case,}\\
&s_n(s,p,\rho,\nu)^{-1}R_{n}\lesssim
\log(n)^{(1-\frac{p}{\rho q})_+} \quad \mbox{in the critical case.}
\end{align*}
\end{theorem}

\begin{theorem} Under condition $C_{up}$ on $Y$, for estimator
  $\hat{f}_n^D$ defined in \eqref{estimateur} and if $q\le p$ in the critical case:
\begin{align*}
\sup_{f \in M(s,p,q,R)} E_{f}(\| \hat
{f}_{n}^D-f \| _{\rho})&\lesssim \big(\frac{\log (n)^{1+\frac{1}{\alpha}}}{n}\big)^{\frac {s} {2s+2\nu +1}}\quad \mbox{in the regular case,}\\
\sup_{f \in M(s,p,q,R)} E_{f}(\| \hat
{f}_{n}^D-f \| _{\rho})&\lesssim \big(\frac{\log (n)^{1+\frac{1}{\alpha}}}{n}\big)^{\frac{s-1/p+1/\rho}{2s+2\nu+1-2/p}} \quad \mbox{in the critical and sparse cases.}\\
\end{align*}
\end{theorem} When the filter satisfies $C_{low}$ and $C_{up}$ the rates of Theorems
$1$ and $2$ match except in the critical case
when $\rho>\frac{p}{q}$, where the upper bound contains an extra logarithmic
factor. This is also observed in density estimation or regression problems (see \cite{Donoho3} and \cite{Donoho2}), and that factor is probably part of the actual rate
of $R_n$: the lower bound is maybe too optimistic.

Analysing the effect of $\nu$, we remark that the rates are similar to the ones established in the white
noise model or other classical non-parametric estimation problems
(examples can be found in \cite{tsybakov03}),
except that here an additional effect reflected by $\nu$ slows the
minimax speed. Indeed the convolution blurs
the observations, making the estimation all the more difficult as $\nu$ is large. This parameter is called ill-posedness
coefficient, explanations about this notion can be found in
\cite{illpose} for example.

Concerning Theorem $3$, we remark that estimator $\hat f_n^D$ is
not optimal first by a log factor in the regular case, which is a
common phenomenon for adaptive estimators as was highlighted in
\cite{tsybakov3}, and secondly by log factors with exponents
proportional to $\frac{1}{\alpha}$. This is due to the difficulty to control
the deviation probability of the estimated wavelet coefficients when
the probability of having small eigenvalues $Y_l$ of the convolution operator is high (ie when $\alpha$ is small).

\bigskip

The main interest of these results is that bounds of the minimax
risk are established in a random operator setting, for a wide scale of $L^{\rho}$
losses, and over general functional spaces which include unsmooth
functions. As far as we know, the lower bound has not been
established in deconvolution problems for such settings even in
the case of deterministic filters.

Let us also note that condition $C_{up}$ imposed on the filter $Y$ is similar to the
conditions generally used in other inverse problems where
the singular values of the
operator are required to decrease polynomially fast. Moreover condition $C_{up}$ concern means of eigenvalues over
diadic blocs, which enables to include filters for
which Fourier coefficients vary erratically individually, but not
in mean, such as some boxcar filters (see \cite{boxcar}). The case of severely ill-posed inverse problems, where the singular
values decrease exponentially fast, has also been studied in
\cite{tsybakov1} for example. 

\section{Simulations}

To illustrate the rates obtained for the upper bound, the
behaviours of estimators $\hat f_n^D$ and $\hat f_n^R$ are examined in practice for the
following settings. We consider the four target functions
(Blocks, Bumps, Heavisine, Doppler) represented on figure $2$, which
were used by Donoho and Johnstone in a series of papers (\cite{dj94} for example). These functions are blurred by
convolution with realizations of a random filter $Y$ and by
adding gaussian noise with root signal to noise ratio ($rsnr$) of three levels: $rsnr \in \{3,5,7\}$. Then the two
estimators are computed in each case and their
performances are examined, judging by the mean square error ($MSE$). For the simulation of the data and the implementation of
the estimators, parts of the WaveD software package
written by Donoho and Raimondo for \cite{deconvolution} were used.

\begin{figure} \centering
\centering
\includegraphics[width=6.5in,height=2in]{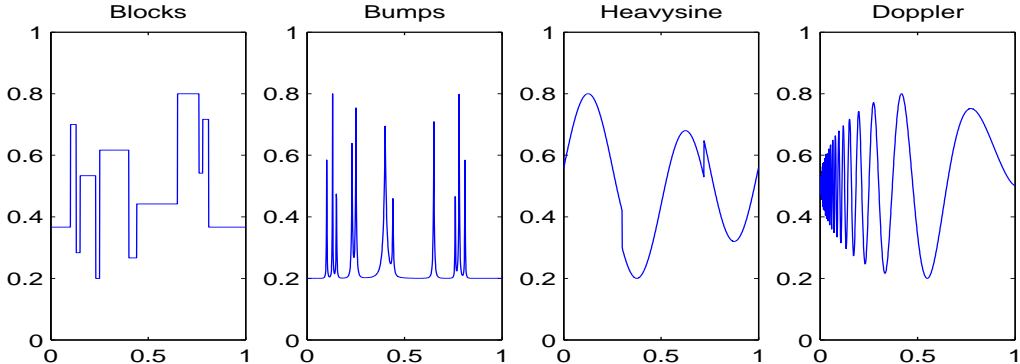}
\caption{\footnotesize Target functions}
\end{figure}

\subsection{Distribution of the filter}

A simple way to represent the blurring effect is the convolution
with a boxcar filter, ie at time $t$ one observes the mean of the
unknown function on an interval $[t-a,t]$ with a random width $a$.
However these kinds of filters have various degrees of ill
posedness depending on $a$. For some numbers called "badly
approximable" numbers, this degree is constant and equal to $3/2$.
For other numbers the situation is more complicated, and the set
of the badly approximable numbers has a Lebesgue measure equal to
zero (more explanations can be found in
\cite{Johnstone} or \cite{deconvolution}). However new results
have been found recently for almost all boxcar widths in \cite{boxcar}
where the near optimal properties of several thresholding
estimators are established.

\bigskip

So as to keep a fixed ill posedness coefficient boxcar
filters are excluded, and one considers convolutions with periodized gamma functions with parameters $\nu$ and $\lambda$:
$$Y(t)= \frac{1}{\int
_{0}^{+\infty}s^{\nu-1}e^{-\lambda s}ds} \sum_{l\in \mathbb{N}} (t+l)^{\nu-1}e^{-\lambda (t+l)},$$ where $\nu$ is a
fixed shape parameter and $\lambda$ is a \textit{random}
scale parameter with a probability distribution function $F_{\alpha}$ parametrized by some $\alpha>0$:

$$F_{\alpha}(t)=\min \bigl(1,2 e^{-\frac {C_{\alpha}} {t^{2\alpha}}}\mathbb{I}(t\ge 0) \bigr),$$ where the constant $C_{\alpha}$ is set such that $E(\lambda)=150$ for all $\alpha$.

\begin{figure} \centering
\centering
\includegraphics[width=6.5in,height=2in]{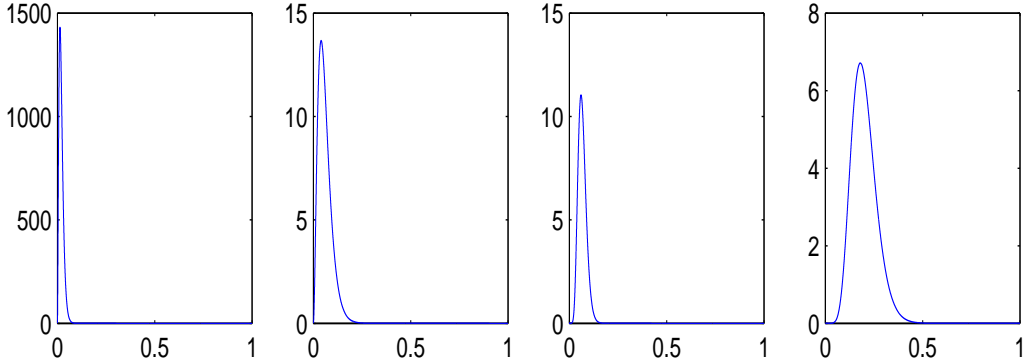}
\caption{\footnotesize Examples of filters, from left to right: $(\nu,\lambda)\in \{(3,150), (3,50), (10,150),(10,50)\}$}
\end{figure}

\bigskip

Such a filter $Y$ satisfies conditions $C_{up}$ and $C_{low}$.
Some examples of its shapes are given in figure $3$: $\nu$ and $\lambda$ can be interpreted respectively as a delay and a spreading parameter. According to the
minimax rates, $f$ should be
(asymptotically) more difficult to estimate for large $\nu$ and for small $\alpha$. This is checked in practice in the next section.

\subsection{Results}

First we focus on the effect of $\nu$ conditionnally to the filter $Y$. An example in medium noise for the
Blocks target is given in figure $4$, where the filter is kept
constant with $\lambda=150$: as expected, both estimators get less and less
efficient when $\nu$ increases. Moreover in practice the thresholds of estimator $\hat f_n^D$ need to be rescaled for each $\nu$, contrarily to
those of estimator $\hat f_n^R$ which is thus more convenient. The same results were obtained
for the other target functions and by examining the $MSE$ of the
estimators, the figures were not included for the sake of conciseness.

\bigskip

Next we set $\nu=1$ and we investigate the effect of the distribution
of the filter $Y$. Both estimators perform well for mean and high
realizations of $\lambda$, but difficulties appear for small
realizations which are all the more frequent as $\alpha$ is small: the
worst case among $10$ simulations is represented in figure $5$ when
$\alpha=2$ and in figure $6$ when $\alpha=0.5$, and the two estimators
perform more poorly in the last case. However they
remain better in that case than a fixed threshold estimator (ie with thresholds completely independent of the filter) also
represented in the figures.

More generally the $MSE$ were computed for several values of $\alpha$
and for the three noise levels.  The results are given in figure $7$:
the shape of the distribution of $Y$ clearly affects estimator
$\hat f_n^D$, and also $\hat f_n^R$ to a much lesser
extent. The smaller $\alpha$, the poorer they behave. Especially the Doppler and Bumps targets are not well
estimated by $\hat f_n^D$ for small $\alpha$, mainly because the high
thresholds make it ignore many of the numerous details of these targets.

\bigskip

Finally estimator $\hat f_n^R$ proves more convenient than estimator $\hat f_n^D$ when the ill-posedness varies, and also less sensitive to the weight of the probability of small eigenvalues.

\begin{figure} [p]\centering
\centering
\includegraphics[width=5in,height=1.5in]{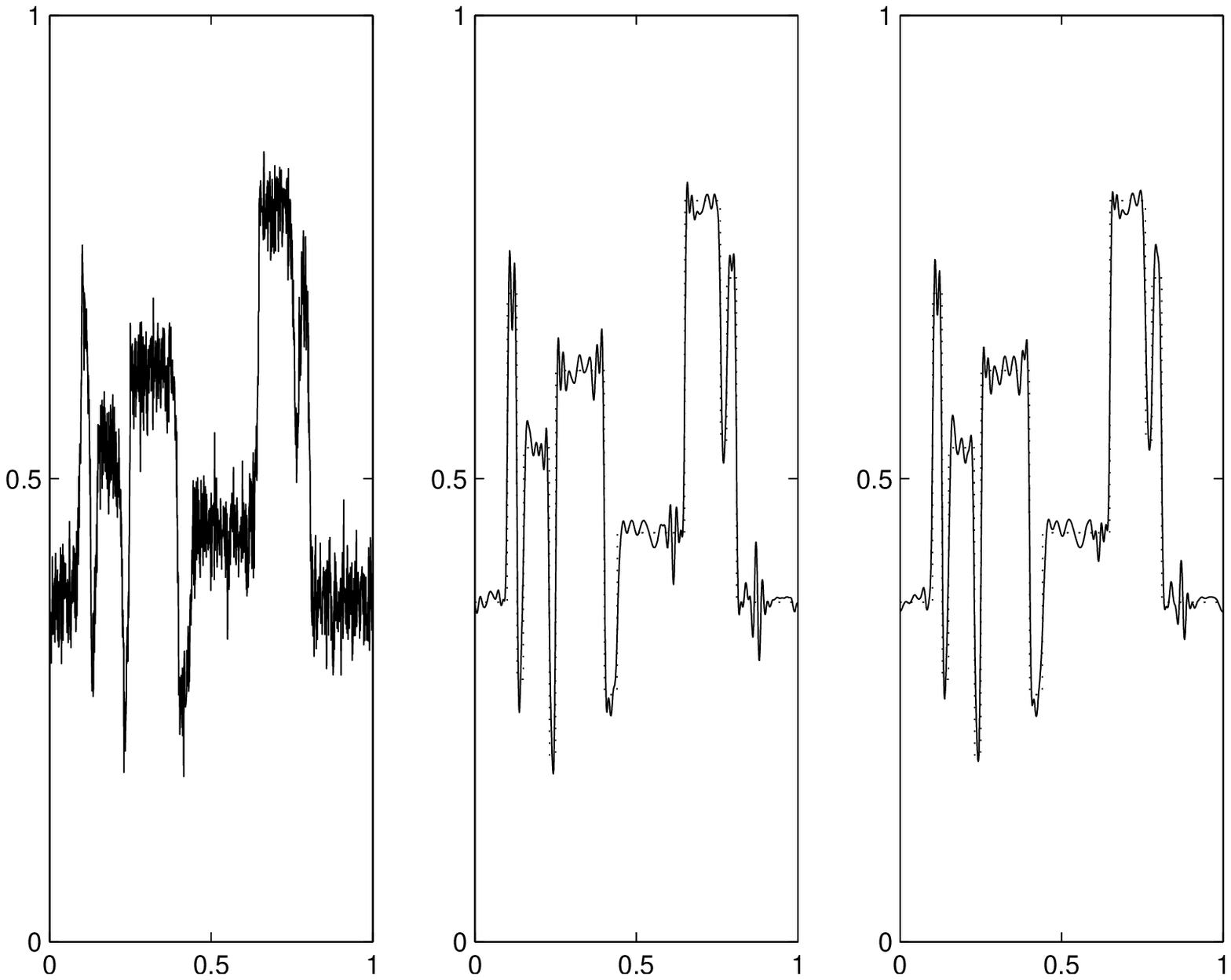}
\includegraphics[width=5in,height=1.5in]{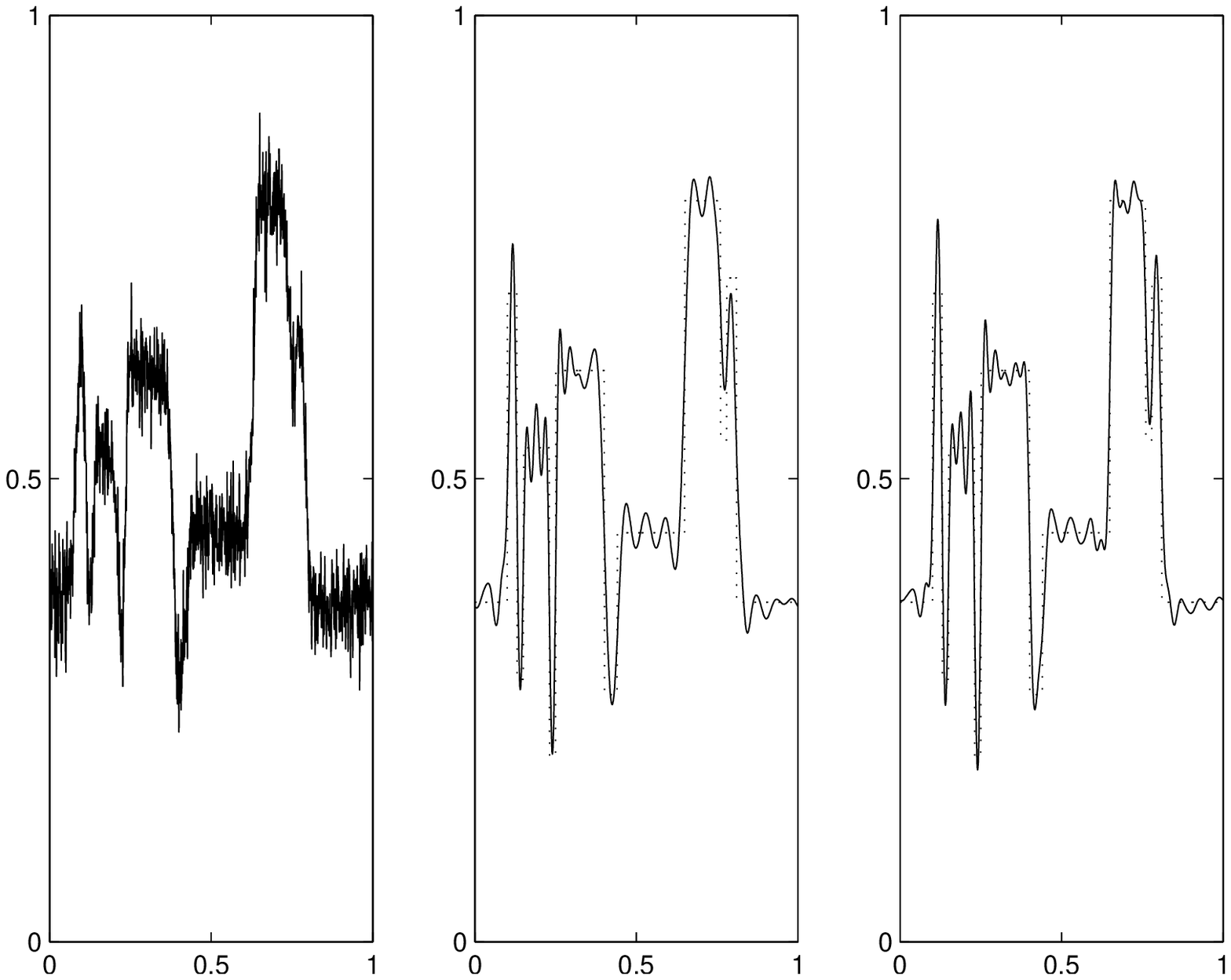}
\includegraphics[width=5in,height=1.5in]{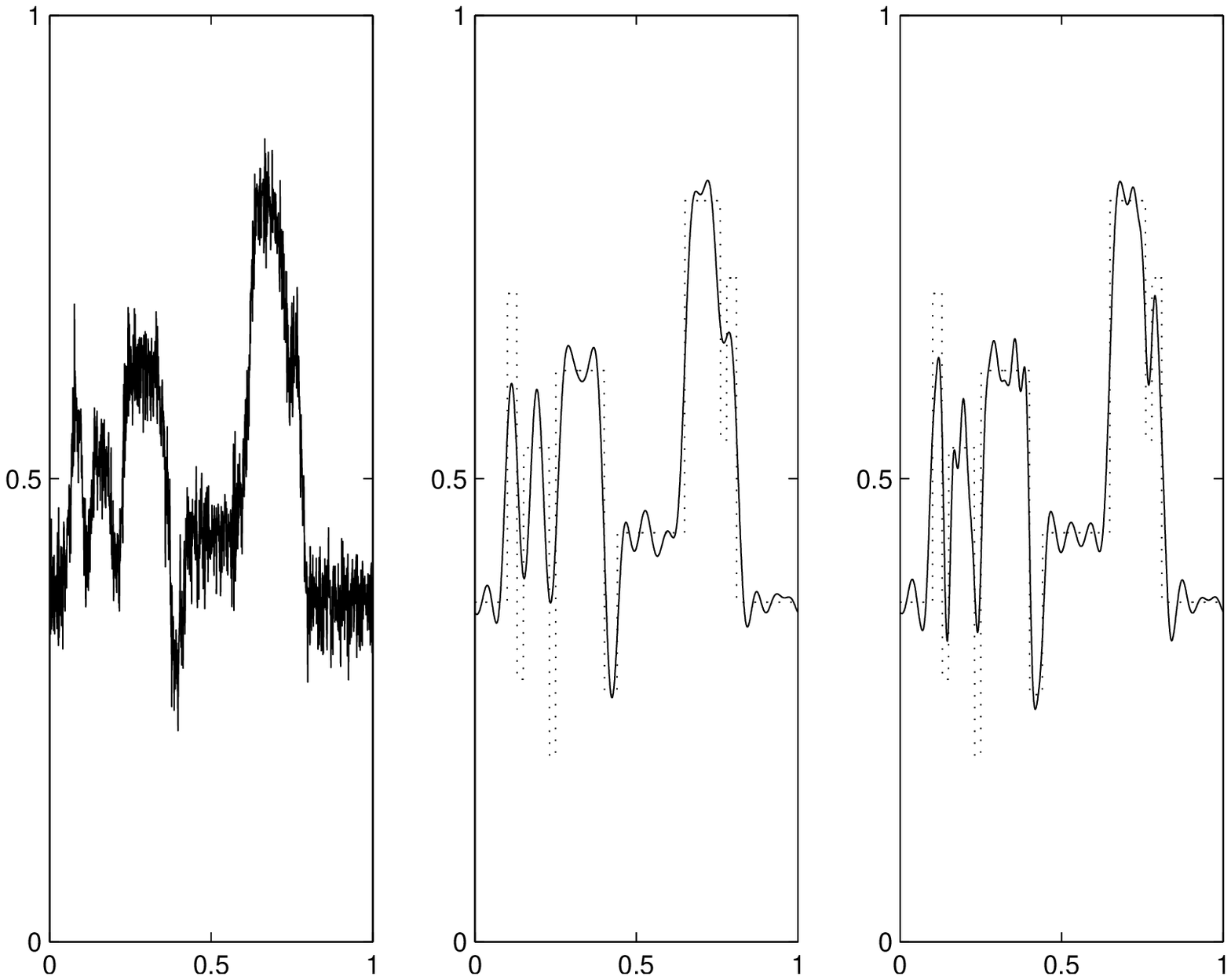}
\caption{\footnotesize Data, estimator $\hat f_n^R$ and estimator $\hat f_n^D$ (left to right) for fixed
  $\lambda=150$ and $\nu=1$ (top), $\nu=3$
  (middle) and $\nu=5$ (bottom)}
\end{figure}

\begin{figure} [p] \centering
\centering
\includegraphics[width=6.5in,height=2in]{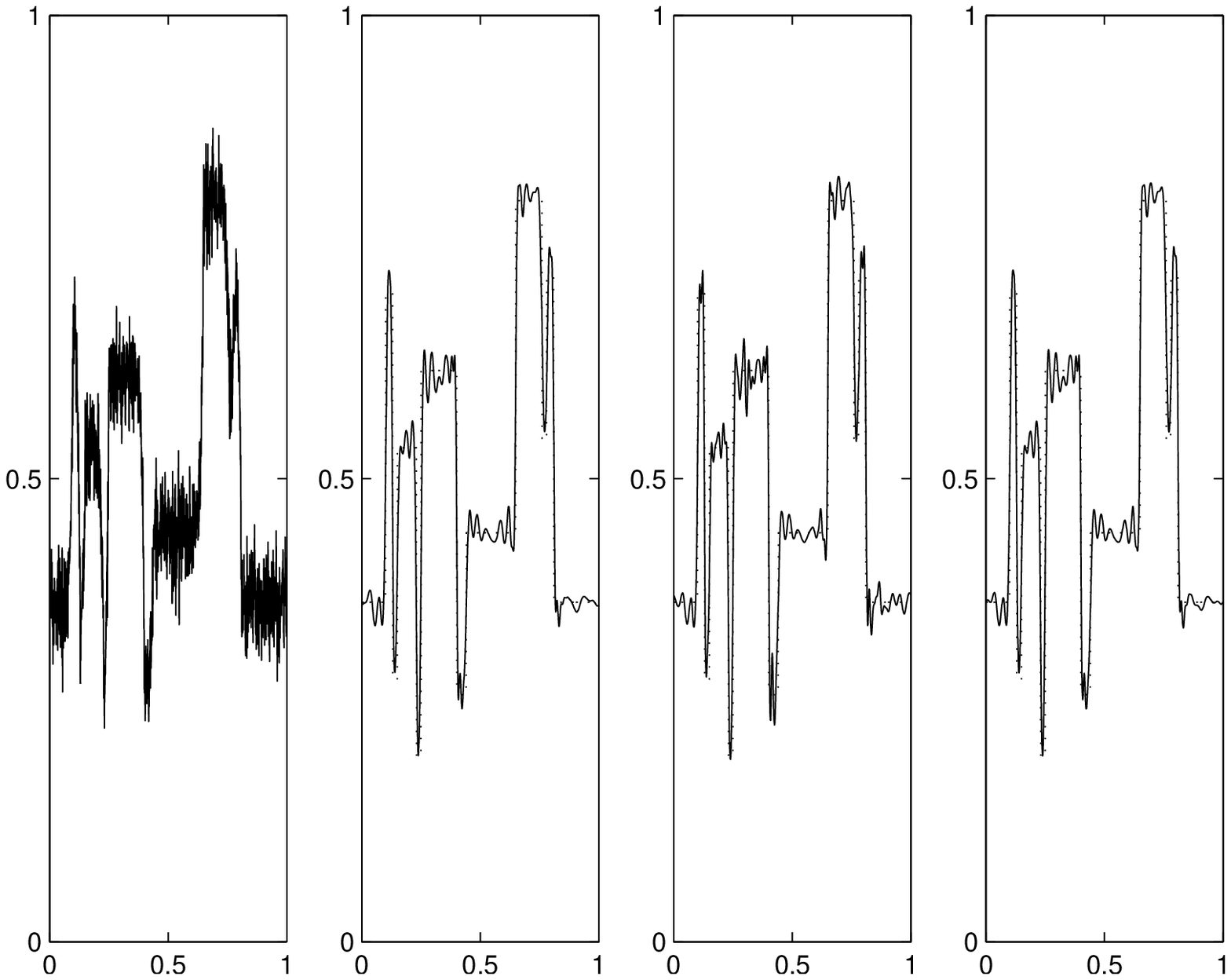}
\caption{\footnotesize Data, estimator $\hat f_n^R$, estimator
  $\hat f_n^D$ and a fixed-threshold estimator (left to right)
  for $\alpha=2$}
\end{figure}

\begin{figure} [p] \centering
\centering
\includegraphics[width=6.5in,height=2in]{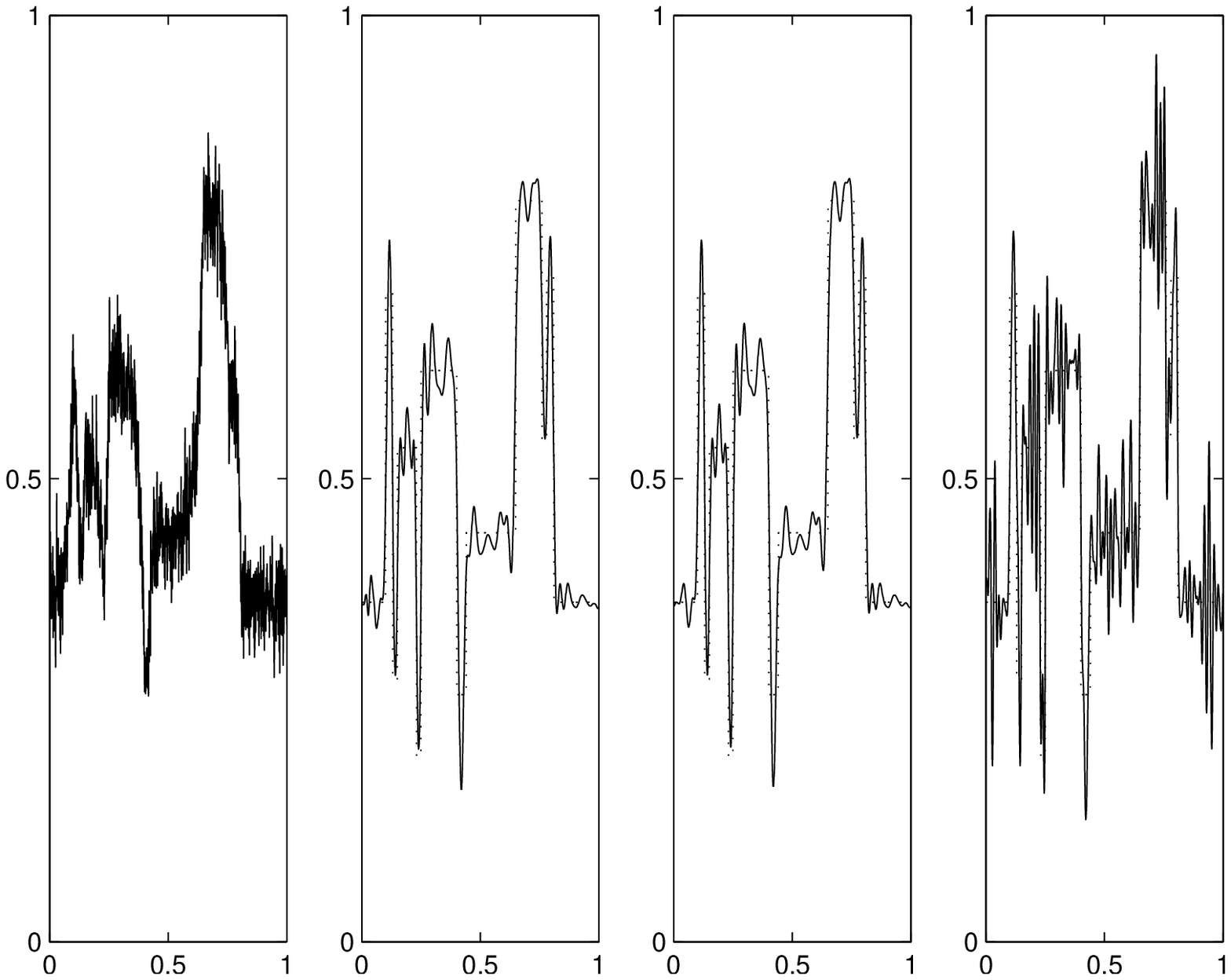}
\caption{\footnotesize Data, estimator $\hat f_n^R$, estimator
  $\hat f_n^D$ and a fixed-threshold estimator (left to right) for $\alpha=0.5$}
\end{figure}

\begin{figure} [p] \centering
\centering
\includegraphics[width=6in,height=6in]{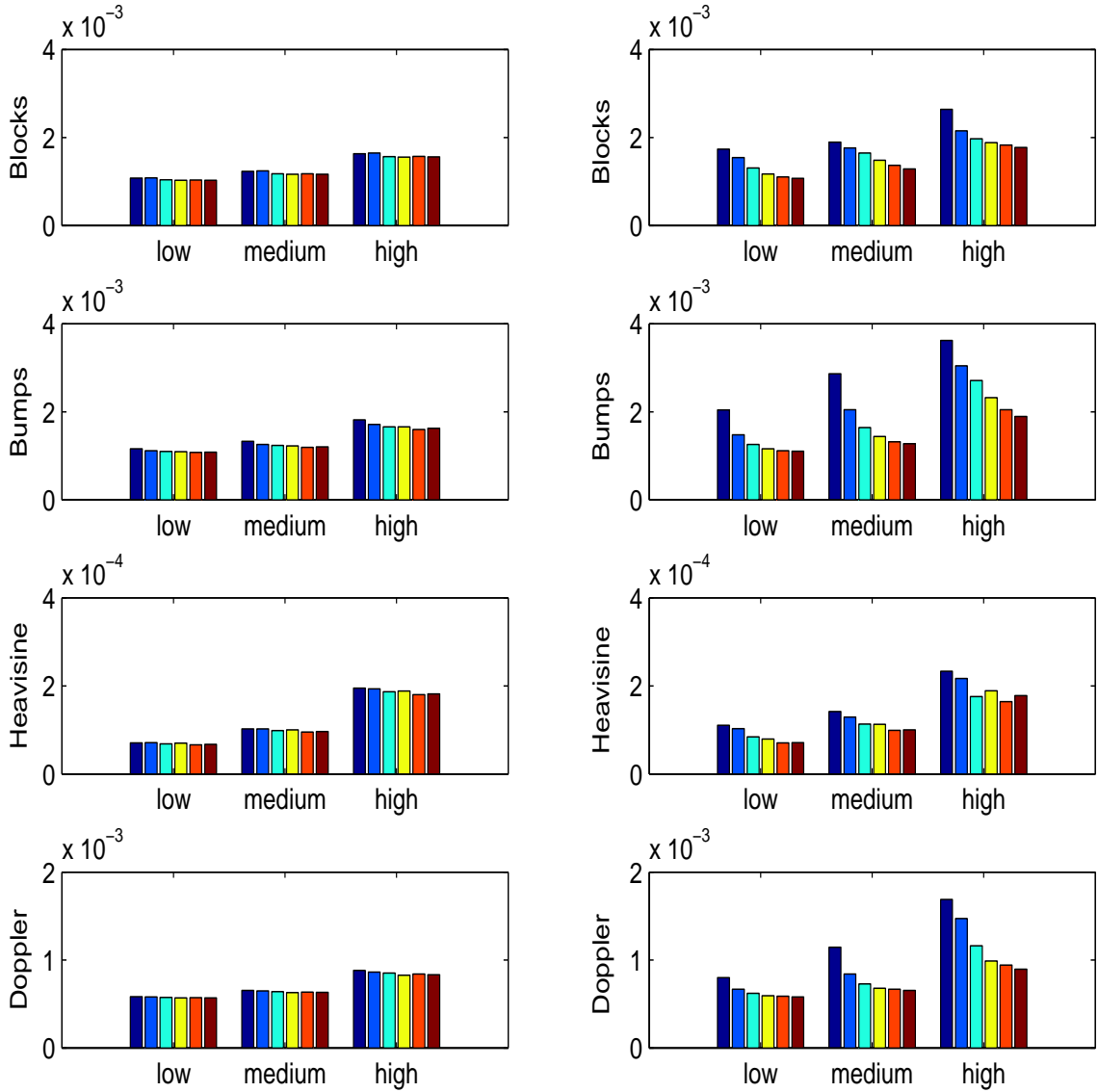}
\caption{\footnotesize Effect of $\alpha$ on the MSE of estimator $\hat f_n^R$ (left) and
  estimator $\hat f_n^D$ (right) for each target, each level of noise and
  $\alpha \in \{0.5,0.6,0.7,0.8,0.9,1 \} $ (left to right in each group) }
\end{figure}

\section{Proofs of the lower and upper bounds}

\subsection{Lower bound}

\subsubsection{Sparse case}

We use a classical lemma on lower bounds (Korostelev and Tsybakov
\cite{korostelev_tsybakov93}):

\begin{lemma}

Let $V$ a functionnal space, $d(.,.)$ a distance on $V$,\\
for $f$, $g$ belonging to $V$ denote by $\Lambda_{n}(f,g)$ the
likelihood ratio : $\Lambda_{n}(f,g)=\frac {dP_{X_{n}^{(f)}}}
{dP_{X_{n}^{(g)}}}$ where $dP_{X_{n}^{(h)}}$ is the probability
distribution of the process $X_{n}$ if $h$ is true.\\
If $V$ contains functions $f_{0},f_{1},\dots , f_{K}$ such that :

\begin{itemize}
\item{$d(f_{k^{'}},f_{k})\geq  \delta >0$ for $k \ne k^{'}$,}
\item{$K\geq \exp (\lambda _{n})$ for some $\lambda _{n}>0$,}
\item{$\Lambda_{n}(f_{0},f_{k})=\exp (z^{k}_{n}-v^{k}_{n})$, where
    $z^{k}_{n}$ is a random variable such that there exists $\pi _{0}>0$ with $P(z^{k}_{n}>0)\geq \pi_{0}$, and $v^{k}_{n}$ are constants,}
\item{$\sup_{k} v^{k}_{n}\leq \lambda _{n}$.}
\end{itemize} Then
$$\sup_{f \in V} P_{X_{n}^{(f)}} \bigl( d(\hat{f}_{n},f)\geq \delta /2
\bigr) \geq \pi _{0}/2,$$ for an arbitrary estimator
$\hat{f}_{n}$.
\end{lemma}


To use this result, we build a finite set of functions belonging
to $M(s,p,q,R)$ as follows. Let $(\psi_{j,k})_{j\ge -1,k\in \mathbb{Z}}$
be an $s-$regular Meyer wavelet basis, which we periodize according to:
$$\Psi_{j,k} (x) =\sum_{l \in \mathbb{Z}}\psi_{j,k} (x+l).$$ In the
sequel we denote by $(\Psi_{j,k})_{(j,k)\in \Lambda}$ the
periodized Meyer wavelet basis obtained this way, where $\Lambda
=\{(j,k)\:|\: j\ge -1; \: k\in R_j\}$ and $R_{j}=\{0,\dots,2^{j}-1\}$.

\bigskip

Now for a fixed level of resolution $j$ set for any $k \in R_{j}$:
$$f_{j,k}=\gamma \Psi _{j,k},$$ with $\gamma \lesssim 2^{-j(s+1/2-1/p)}$ such that $\| f_{j,k} \| _{s,p,q} \leq R
$. Set also $f_{0}=0$.

\bigskip

Let us choose for $d$ the distance $d(f,g)=\| f-g \|_{\rho}$.
Because of the relation between the $L^{\rho}$ norm of a linear
combination of wavelets of fixed resolution $j$ and the $l^{\rho}$
norm of the corresponding coefficients (see \cite{Meyer}), we have for any $k,k^{'}
\in R_{j}$, $k \ne k^{'}$:
$$d(f_{j,k^{'}},f_{j,k})=\| \gamma \Psi _{j,k^{'}}- \gamma \Psi _{j,k}
\|_{L^{\rho }}\asymp \gamma 2^{j(1/2-1/\rho)}.$$ In this framework
we have : $K=2^{j}$ and $\delta\asymp \gamma 2^{j(1/2-1/\rho)}$.
So as to apply the lemma, we have to find parameters $\gamma (n)$
and $j(n)$ such that the other hypotheses of the lemma are
satisfied, which will be true if :
$$P_{f_{j,k}} \Bigl( \ln (\Lambda_{n}(f_{0},f_{j,k}))
\geq -j(n)\ln(2) \Bigr) \geq \pi _{0}> 0,$$ uniformly for all
$f_{j,k}$. Moreover we have :
\begin{align*}
P_{f_{j,k}} \Bigl( \ln (\Lambda_{n}(f_{0},f_{j,k})) \geq
-j(n)\ln(2) \Bigr) 
&\ge1-P_{f_{j,k}} \Bigl( |\ln (\Lambda_{n}(f_{0},f_{j,k}))| > j(n)\ln(2) \Bigr)\\
&\geq 1-E_{f_{j,k}} \Bigl( |ln(\Lambda_{n}(f_{0},f_{j,k})| \Bigr)
/(j(n)\ln(2)).
\end{align*} So the previous condition is satisfied when $\gamma (n)$ and
$j(n)$ are chosen such that, with a constant $0<c<1$:
\begin{equation}
\label{cond} E_{f_{j,k}}\Bigl( |ln(\Lambda_{n}(f_{0},f_{j,k}))|
\Bigr) \leq cj(n)\ln(2).
\end{equation}

\bigskip

Consider two
hypotheses $f_{0}$ and $f_{j,k}$, and let us determine the
likelihood ratio of the corresponding distributions of the observations
$(X_{n}(t),Y(t))_{t\in[0,1]}$. Let $F$ be a bounded measurable
function. Since $Y$ is assumed to be independent of $W$ and free with
respect to $f$
in \eqref{eqdep}, we have:
\begin{align*}
E_{f_{j,k}}\big[F\big(X^n,Y\big)\big]&=E\big[E\{F\big( (\int _{0}^{t}f_{j,k}\star Y(s)ds+\sigma
n^{-1/2}W(t),Y(t))_{t\in[0,1]}\big)\quad | Y \} \big]\\
&=\int E\{F\big( \sigma n^{-1/2}\tilde{W},y\big)\}dP_{Y}(y),
\end{align*} where $P_{Y}$ denotes the distribution of $Y$ and $\tilde{W}(t)=W(t)+\int
_{0}^{t}\sigma^{-1}n^{1/2}f_{j,k}\star y(s)ds.$

\bigskip

For a given function $y$ let $h_{j,k}^{y}$ be defined by:
$h_{j,k}^{y}(t)=\sigma^{-1} n^{1/2} f_{j,k} \star y (t).$ We
assumed that $Y$ takes its values in $L^{2}([0,1])$ so for each of
its realization there exists a constant $C_y$ such that for all
$t\in [0,1]$, $\int_{0}^{t}(h_{j,k}^{y})^{2}(s)ds < C_y $ and we can apply
the formula of Girsanov: the process $\tilde{W}$ is a Wiener process under the probability
$Q$ defined by $$dQ=\exp \big[-\int_{0}^{1}
h_{j,k}^{y}(t)dW(t)-\frac{1}{2}\int_{0}^{1} (h_{j,k}^{y}(t))^{2}dt\big]
dP.$$ Thus for any function $y$:
\begin{align*}
E_{P}\big[F\big( \sigma
n^{-1/2}\tilde{W},y\big)\big]&=E_{Q}\big[F\big( \sigma
n^{-1/2}\tilde{W},y\big) \exp \big[\int_{0}^{1}
h_{j,k}^{y}(t)dW(t)+\frac{1}{2}\int_{0}^{1} (h_{j,k}^{y}(t))^{2}dt\big]
\big]\\&=E_{Q}\big[F\big( \sigma n^{-1/2}\tilde{W},y\big) \exp
\big[\int_{0}^{1} h_{j,k}^{y}(t)d\tilde{W}(t)-\frac{1}{2}\int_{0}^{1}
(h_{j,k}^{y}(t))^{2}dt\big] \big]\\
&=E_{P}\big[F\big( \sigma n^{-1/2}W,y\big) \exp \big[\int_{0}^{1}
h_{j,k}^{y}(t)dW(t)-\frac{1}{2}\int_{0}^{1} (h_{j,k}^{y}(t))^{2}dt\big]
\big].
\end{align*} So finally:
$$\Lambda_{n}(f_{0},f_{j,k})=\exp \big[-\int_{0}^{1}
\frac{f_{j,k} \star Y (t)}{\sigma
n^{-1/2}}dW(t)+\frac{1}{2}\int_{0}^{1} \big(\frac{f_{j,k} \star Y
(t)}{\sigma n^{-1/2}}\big)^{2}dt\big].$$

\bigskip

We can now examine under which conditions \eqref{cond} is true. We
have:
\begin{align*}
&E|\ln (\Lambda_{n}(f_{0},f_{j,k}))|=E|\frac{\gamma n^{1/2}}{\sigma
}\int_0 ^1 \Psi_{j,k} \star Y
(t)dW(t)-\frac{\gamma^{2}n}{2\sigma^{2}}\int_0 ^1 (\Psi_{j,k} \star Y
(t))^{2}dt|\leq A_{n} + B_{n}, \mbox{  with:}\\
&B_{n}=\frac{\gamma^{2}
n}{2\sigma^{2}}E\big(\int_0 ^1 (\Psi_{j,k} \star Y
(t))^{2}dt \big),\\
&A_{n}=\frac{\gamma n^{1/2}}{\sigma}E|\int_0 ^1\Psi_{j,k} \star Y (t) dW(t)|\leq\frac{\gamma n^{1/2}}{\sigma}\big(E(\int_0 ^1
\Psi_{j,k} \star Y (t) dW(t))^{2} \big)^{1/2}
\leq (2B_{n})^{1/2},
\end{align*} where we used Jensen's inequality for $A_n$.

\bigskip

Let us find a bound for $B_n$. We introduce the Fourier
coefficients of $Y$ and $\Psi_{j,k}$ denoted by $Y_{l}$ and
$\Psi_{j,k,l}$ for all $l\in \mathbb{Z}$. Since the Fourier Transform
of $\Psi_{j,k}$ is bounded by $2^{-j/2}$ we have:
$$B_{n}=\frac{\gamma^{2}n}{2\sigma^{2}}E_{f_{j,k}}\big(\frac{1}{2\pi} \sum_{l\in
\mathbb{Z}} |Y_{l}\Psi_{j,k,l}|^{2}\big) \lesssim \gamma^{2}n
2^{-j} E_{f_{j,k}}\big(\sum_{l\in C_{j}} |Y_{l}|^{2}\big),$$ where
$C_{j}$ is the set of integers where the coefficients
$\Psi_{j,k,l}$ are not equal to zero (it can easily be shown that this
set does not depend on $k$).

\bigskip

The support of the Fourier transform of the Meyer wavelet is
included in $[-\frac{2 \pi}{3},-\frac{8 \pi}{3}]\cup[\frac{2
\pi}{3},\frac{8 \pi}{3}]$. So $\Psi_{j,k,l}=0$ as soon as $|2 \pi
2^{-j}l| \in [\frac{2 \pi}{3},\frac{8 \pi}{3}]^{c}$, and
$C_{j}\subset [-2^{j+1},-2^{j-2}]\cup[2^{j-2},2^{j+1}]$ for all
$j$. Then under condition $C_{low}$ and noticing that $Y_{-l}=Y_l$ we obtain:

$$B_{n} \lesssim \gamma^{2}n 2^{-2\nu j}.$$ Finally, condition \eqref{cond} holds if we choose $\gamma$
and $j$ such that:

$$\gamma^{2}n 2^{-2\nu j} \lesssim j,\quad \mbox{ and } \quad \gamma \lesssim 2^{-j(s+1/2-1/p)}.$$ We choose the following values that satisfy those two conditions:
$$\gamma \asymp 2^{-j(s+1/2-1/p)}, \quad \mbox{ and } \quad  2^{j} \asymp (n/log(n))^{1/(2s+2\nu +1-2/p)}.$$ Finally, using the lemma and the inequality of Markov, for $\sigma ((X_{n}(t),Y(t)),\:
t\in[0,1])-$measurable estimators $\hat {f}_{n}$ the following
bound holds:

$$\inf_{\hat {f}_{n}} \sup_{f \in M(s,p,q,S)} E_{f}(\| \hat {f_{n}}-f \|
_{\rho}) \gtrsim \gamma 2^{j(1/2-1/\rho)} \asymp \big(\frac{\log(n)}{n}\big)^{\frac{s-1/p+1/\rho}{2s+2\nu +1-2/p}}.$$

\subsubsection{Regular case}

Here we consider another set of functions belonging to
$M(s,p,q,R)$. We use the periodized Meyer wavelet basis $(\Psi
_{j,k})$ like before. But now we set for any $\epsilon \in \{
-1,+1 \}^{R_{j}}$:
$$f_{j,\epsilon }=\gamma \sum_{k \in R_{j}}\epsilon _{k}\Psi
_{j,k},$$ with $\gamma \lesssim 2^{-j(s+1/2)}$ such that $\|
f_{j,\epsilon } \| _{s,p,q} \leq S $. We also set
$I_{j,k}=[\frac{k}{2^{j}},\frac{k+1}{2^{j}}].$

\bigskip

We use an adaptation of lemma $10.2$ in \cite{wavelets} to the
case of Meyer wavelets (that do not have compact supports) and of the
norm $\|.\|_{\rho}$:

\begin{lemma}
Suppose the likelihood ratio satisfies for some constant $\lambda
$:
$$P_{f_{j,\epsilon }}\big(\Lambda_{n}(f_{j,\epsilon ^{k}},f_{j,
\epsilon })\geq e^{-\lambda}\big)\geq p_{*}> 0,$$ uniformly for
all $f_{j,\epsilon }$ and all $k \in R_{j}$, where $\epsilon ^{k}$
is equal to $\epsilon$ except for the $k^{th}$ element which is
multiplied by $-1$. Then the following bound holds:
$$\max_{\epsilon \in \{ -1,+1 \} ^{R_{j}}} E_{f_{j,\epsilon }}(\| \hat
{f_{n}} - f_{j,\epsilon } \| _{\rho}) \geq C
2^{j/2}\gamma e^{-\lambda}p_{*},$$ where $C$ is
positive and depends only on $\rho$.
\end{lemma}

\bigskip

Similarly to the sparse case, the hypothesis of this lemma is
satisfied if, for a small enough constant $c$:
$$E_{f_{j,\epsilon }}|\ln \big(\Lambda_{n}(f_{j,\epsilon ^{k}},f_{j,
\epsilon }) \big)|\le c.$$

\bigskip

\noindent Now the log-likelihood is equal to:

$$\ln \big(\Lambda_{n}(f_{j,\epsilon ^{k}},f_{j,\epsilon
})\big)=\frac{2\gamma n^{1/2}}{\sigma}\int_0 ^1\Psi_{j,k} \star Y
(t)dW(t)-\frac{2\gamma^{2}n}{\sigma^{2}}\int_0 ^1 [\Psi_{j,k}
\star Y (t)]^{2}dt.$$ Like before, we only need to dominate the
following quantity:

$$B_{n}=\gamma^{2}nE_{f_{j,\epsilon }}(\int_0 ^1
(\Psi_{j,k} \star Y (t))^{2}dt).$$

\noindent We use the same bound as in the sparse case, under assumption
$C_{low}$. The parameters have to be chosen such
that:
$$\gamma^{2} n 2^{-2\nu j} \lesssim 1 \quad \mbox{ and } \quad \gamma
\lesssim 2^{-j(s+1/2)}.$$ Finally the regular rate is obtained for the
following choices:
$$\gamma \asymp 2^{-j(s+1/2)}, \quad \mbox{ and } \quad 2^{j} \asymp
n^{1/(2s+2\nu +1)}.$$

\bigskip

\begin{proof} \textit{of the lemma}\\

\bigskip

The Meyer wavelet satisfies $\exists A>0$ such that $|\psi(x)|\le \frac{A}{1+|x|^{2}}$. Consequently:
\begin{align*}
\bigl(\int_{I_{j,k}}|\Psi_{j,k}(x)dx|^{\rho}
\bigr)^{1/\rho}&=2^{j(\frac{1}{2}-\frac{1}{\rho})}\big(\int_{0}^{1}|\sum_{l
\in \mathbb{Z}}\psi (x+2^{j}l)|^{\rho}dx \big)^{1/\rho}\\
&\ge 2^{j(\frac{1}{2}-\frac{1}{\rho})}\big(
\int_{0}^{1}|\psi(x)|^{\rho}dx-\sum_{l \in
\mathbb{Z}^{*}}\int_{0}^{1}|\psi(x+2^{j}l)|^{\rho}dx \big)^{1/\rho}\\
&\ge 2^{j(\frac{1}{2}-\frac{1}{\rho})}\big(
\int_{0}^{1}|\psi(x)|^{\rho}dx-\frac{ A^{\rho}}{2^{2\rho
j}}\sum_{l \in \mathbb{N}^{*}}\frac{1}{(l/2)^{2\rho}} \big)^{1/\rho}\\
&\ge c 2^{j(\frac{1}{2}-\frac{1}{\rho})},
\end{align*} for $j$ large enough and $c>0$ depends only on $\rho$.

\bigskip

Then using a concavity inequality and similar arguments as in the compact support case, we have:
\begin{align*}
&\max _{\epsilon} E_{f_{j,\epsilon }}(\| \hat {f_{n}} -
f_{j,\epsilon } \| _{\rho}) \ge 2^{-2^{j}}
\sum_{\epsilon}
E_{f_{j,\epsilon}}[\sum_{k=0}^{2^{j}-1} \int_{I_{j,k}}|\hat {f_{n}} - f_{j,\epsilon
}|^{\rho}]^{\frac{1}{\rho}}\\
&\ge 2^{-2^{j}+j(\frac{1}{\rho}-1)} \sum_{\epsilon}
\sum_{k=0}^{2^{j}-1} E_{f_{j,\epsilon}}[\int_{I_{j,k}}|\hat {f_{n}} - f_{j,\epsilon
}|^{\rho}]^{\frac{1}{\rho}}\\
&\ge 2^{-2^{j}+j(\frac{1}{\rho}-1)} \sum_{k=0}^{2^{j}-1} \sum_{\epsilon | \epsilon_{k}=1}
E_{f_{j,\epsilon}}[(\int_{I_{j,k}}|\hat {f_{n}} - f_{j,\epsilon
}|^{\rho })^{\frac{1}{\rho}}+ \Lambda_{n}(f_{j,\epsilon ^{k}},f_{j, \epsilon })
(\int_{I_{j,k}}|\hat {f_{n}} - f_{j,\epsilon ^{k}}|^{\rho
})^{\frac{1}{\rho}}]\\
&\ge 2^{-2^{j}+j(\frac{1}{\rho}-1)} \sum_{k=0}^{2^{j}-1} \sum_{\epsilon |
\epsilon_{k}=1} E_{f_{j,\epsilon}}[\delta I\{
\int_{I_{j,k}}|\hat {f_{n}} - f_{j,\epsilon }|^{\rho } \ge
\delta^{\rho} \}+ \Lambda_{n}(f_{j,\epsilon ^{k}},f_{j, \epsilon
})\delta I\{\int_{I_{j,k}}|\hat {f_{n}} - f_{j,\epsilon^{k}
}|^{\rho } \ge \delta^{\rho} \}]
\end{align*} with $\delta=c \gamma
2^{j(\frac{1}{2}-\frac{1}{\rho})}.$

Noticing that $$\big(\int_{I_{j,k}}|\hat {f_{n}} - f_{j,\epsilon
}|^{\rho }\big)^{1/\rho}+\big(\int_{I_{j,k}}|\hat {f_{n}} -
f_{j,\epsilon ^{k}}|^{\rho }\big)^{1/\rho} \ge 2 \gamma
\big(\int_{I_{j,k}}|\Psi_{j,k}(x)|^{\rho}\big)^{1/\rho} \ge 2
\gamma c 2^{j(\frac{1}{2}-\frac{1}{\rho})}$$ for $j$ large enough,
the end of the proof follows as in \cite{wavelets}.
\end{proof}

\subsection{Upper bounds}

\subsubsection{Properties of the estimated wavelet coefficients}

The performances of the thresholding estimators rest on the properties of the estimated wavelet coefficients
$\hat{\beta}_{j,k}$. In the sequel we will also need properties for the
estimators $\hat{\alpha}_{j,k}$ defined the same way as
$\hat{\beta}_{j,k}$ in estimator \eqref{estimateur} except with
$\Phi$ instead of $\Psi$. We have the following results:

\begin{proposition} Under condition $C_{up}$ we have for all $j\ge
  -1$, $k \in R_j$ and $r>0$,
$$E(|\hat{\beta}_{j,k}-\beta_{j,k}|^{r})\lesssim \bigr(\frac{2^{\nu
    j}}{\sqrt{n}}\bigl)^{r} \quad \mbox{and} \quad E(|\hat{\alpha}_{j,k}-\alpha_{j,k}|^{r})\lesssim \bigr(\frac{2^{\nu
    j}}{\sqrt{n}}\bigl)^{r},$$ and there exist positive constants
    $\kappa,$ and $\kappa'$
    such that for all $\lambda \ge 1$,
$$P(|\hat{\beta}_{j,k}-\beta_{j,k}|\ge \frac{2^{\nu
    j}}{\sqrt{n}}\lambda) \lesssim 2^{-\kappa \lambda^{\frac{2\alpha}{\alpha+1}}} \quad \mbox{and} \quad P(|\hat{\beta}_{j,k}-\beta_{j,k}|\ge \sqrt{\frac{U_j^Y}{n}}\lambda) \lesssim 2^{-\kappa' \lambda^2},$$  where the constants in the inequalities do not
depend on $j$, $k$ and $\lambda$.
\end{proposition}

\bigskip

\begin{proof} \textit{of Proposition $1$}\\

\bigskip

Remark that conditionally to the process $Y$,
$(\hat{\beta}_{j,k}-\beta_{j,k})$ is a centered gaussian variable
with variance:
$$Var (|\hat{\beta}_{j,k}-\beta_{j,k}| \quad |Y)=E[\frac{\sigma^2}{n}\sum _{l
  \in \mathbb{Z}} |\frac {W_l}{Y_l} \Psi_{j,k,l}|^{2}\quad |Y].$$
Since the Fourier transform of the Meyer wavelet is bounded by
$2^{-j/2}$ and only\\ $l\in [-(2^{j+1}-1),-2^{j-2}]\cup
[2^{j-2},2^{j+1}-1]$ has to be considered, we
have for some constant $C>0$:
$$Var (|\hat{\beta}_{j,k}-\beta_{j,k}| \quad |Y)\le C
U_j^{Y}/n.$$

\bigskip

Thus the moment of order $r$ of
$(\hat{\beta}_{j,k}-\beta_{j,k})$ is bounded by
$$E(|\hat{\beta}_{j,k}-\beta_{j,k}|^{r})\lesssim E[(Var (|\hat{\beta}_{j,k}-\beta_{j,k}| \quad |Y))^{r/2}]
\lesssim E[(U_{j}^{Y}/n)^{r/2}],$$ and by similar arguments the same bound
holds for $(\hat{\alpha}_{j,k}-\alpha_{j,k})$ because the support of the
Fourier Transform of $\phi_{j,k}$ is $\frac{4\pi}{3}[-2^j, 2^j]$.

\bigskip

For the deviation probability we use a probabilistic inequality for a
centered standard gaussian variable $Z$. Conditionally to $Y$ we have:
\begin{align*}
P(|\hat{\beta}_{j,k}-\beta_{j,k}|> \frac{2^{\nu
    j}}{\sqrt{n}}\lambda \quad |Y)
&\le P(|Z|\ge \lambda  \sqrt{2^{2\nu j} /(CU_{j}^{Y})} \quad |Y)\\
&\lesssim \frac{1}{\lambda  \sqrt{2^{2\nu j} /(CU_{j}^{Y})}}
\exp (-\frac{\lambda^2 2^{2\nu j}}{2CU_{j}^{Y}}).
\end{align*} Then we take the expectation over $Y$, by Cauchy Schwartz
    we obtain for $\lambda\ge 1$:
\begin{align*}
P(|\hat{\beta}_{j,k}-\beta_{j,k}|> \frac{2^{\nu
    j}}{\sqrt{n}}\lambda) &\lesssim
\sqrt{ E(\frac{U_j^{Y}}{2^{2\nu j}})
E(\exp (-\frac{\lambda^2 2^{2\nu j}}{CU_{j}^{Y}}})).
\end{align*}

The end of the proof is directly deducible from the lemma
below, and the last part of Proposition $1$ is easily proved by
replacing $2^{\nu j}$ by $\sqrt{U_{j}^{Y}}$ in the three inequalities above.

\begin{lemma}
Let $X_j$ be the following random variable: $X_j=\frac{U_{j}^{Y}}{2^{2
\nu j}}$. For all $j\ge 0$ there exists positive constants $C'$,
$C''$, $C(.)$ such that for all $r>0$:
\begin{align*}
E(e^{-\frac{r}{X_{j}}})\le C' e^{-C''r^{\frac{\alpha}
{\alpha+1}} },\quad \mbox{ and } \quad
E(X_{j}^{r})\le C(r).
\end{align*}
\end{lemma}

\bigskip

\begin{proof} \textit{of the lemma}\\

\bigskip

\noindent For all $r>0$ we have:
\begin{align*}
E(e^{-\frac{r}{X_{j}}})&=\int
_{0}^{1}P(e^{-\frac{r}{X_{j}}}\ge u)du\\
&=r \int
_{0}^{+\infty}P(X_{j}\ge 1/u)e^{-r u}du\\
&\le r \int _{0}^{1}P(X_{j}\ge 1/u)e^{-r
u}du +e^{-r}\\
&\lesssim r \int _{0}^{1}
e^{-r u -c/u^{\alpha} }du +e^{-r},
\end{align*} and one can check that there exists $C''>0$ such that  $\int _{0}^{1}
e^{-r u -c/u^{\alpha} }du \lesssim e^{-C''r^{\frac{\alpha} {\alpha+1}} }$.

\bigskip

\noindent The second part of the lemma is easily proved by using
similar arguments.

\end{proof}
\end{proof}

\subsubsection{Proof of the sharp rates}

In the regular and critical zones, estimator \eqref{estimateur} is not
optimal up to a logarithmic factor. In order to show that
the rates of Theorem $1$ are sharp, we exhibit estimators achieving the
rates of Theorem $2$. Those are not as interesting in practice as \eqref{estimateur}, since they depend
on caracteristics of $f$, ie they are not adaptive.

\bigskip

We will use the following bound to estimate the risks, which holds for any $-1 \le
j_m \le j_M \le \infty$ and any set of random or deterministic coefficients
$\tilde{\beta}_{j,k}$ such that the quantities below are finite:

\begin{equation}
\label{norme}
E\| \sum_{j_m\le j \le j_M} \sum_{k \in R_j} \tilde{\beta}_{j,k} \Psi_{j,k} \|_{\rho} \lesssim
\sum_{j_m\le j \le j_M} 2^{j(\frac{1}{2}-\frac{1}{\rho})}\bigl(\sum_{k
  \in R_j} E| \tilde{\beta}_{j,k}|^{\rho} \bigr)^{\frac{1}{\rho}}.
\end{equation} The proof is immediate by Minkowski inequality, the
fact that $\| \sum_{k \in R_j} \tilde{\beta}_{j,k} \Psi_{j,k} \|_{\rho}
\asymp  2^{j(\frac{1}{2}-\frac{1}{\rho})} \|\tilde{\beta}_{j,.}\|_{l_{\rho}}$ (established in \cite{Meyer}) and a concavity argument.

\bigskip

Let us denote: $\nu '= \nu+1/2$ and $\epsilon = ps- \nu'(\rho-p)$. We distinguish
two cases: $\rho \le p$ and $p<\rho$. In the first case $M(s,p,q,R)$ is included in the regular zone. By
concavity we have: $$\inf_{\hat {f}_{n}} \sup_{f \in M(s,p,q,R)} E_{f}\| \hat
{f}_{n}-f \| _{\rho} \le
\inf_{\hat {f}_{n}} \sup_{f \in M(s,p,q,R)} E_{f}\| \hat
{f}_{n}-f \| _{p}.$$ So seeing the expected rate only the case $\rho=p$
needs to be considered. We take the following linear
estimator:
\begin{equation*}
\hat f _n=\sum_{k \in R_{j_1}} \hat{\alpha}_{j_1,k} \Phi_{j_1,k}.
\end{equation*}

For any $f\in M(s,p,q,R)$ the risk is composed of a bias error and a
stochastic error:
$$E_{f}\| \hat{f}_{n}-f \| _{p} \le A_s +A_s,$$ with:
\begin{align*}
A_s=E\|\sum_{k \in R_{{j_1}}}
(\hat{\alpha}_{j_1,k}-\alpha_{j_{1},k}) \Phi_{j_1,k}\|_{p} \lesssim 2^{j_1
  (\frac{1}{2}-\frac{1}{p})} [\sum _{k \in R_{j_{1}}}
E|\hat{\alpha}_{j_1,k}-\alpha_{j_1,k}|^{p}]^{\frac{1}{p}}\lesssim \bigr(\frac{2^{\nu j_1}}{ \sqrt{n}} \bigl)2^{
  \frac{j_1}{2}}=\frac{2^{\nu' j_1}}{ \sqrt{n}},
\end{align*}

\begin{align*}
A_b=\| \sum _{j >
j_1}\sum _{k \in R_j} \beta_{j,k} \Psi_{j,k}\|_{p} \lesssim \sum
_{j>j_1} 2^{j (\frac{1}{2}-\frac{1}{p})} (\sum _{k \in R_j}
|\beta_{j,k}|^{p})^{\frac{1}{p}}\lesssim \sum
_{j>j_1} 2^{j (\frac{1}{2}-\frac{1}{p})}2^{-j
(s+\frac{1}{2}-\frac{1}{p})}\lesssim 2^{-j_1 s},
\end{align*} and we obtain the rate by choosing
$j_1=[\frac{log_2 (n)}{2s+2 \nu '}]$.

\bigskip

In the second case ($p<\rho$) we consider the following estimator:

$$\hat f _n=\sum_{k \in R_{j_1+1}} \hat{\alpha}_{j_1+1,k} \Phi_{j_1+1,k} +\sum_{j_1<j<j_2} \sum_{k \in R_j} \hat{\beta}_{j,k} I_{ \{|\hat{\beta}_{j,k}| \geq \lambda _j\}} \Psi_{j,k},$$
 where:
$$2^{j_1}\approx n^{\frac{1}{2s+2 \nu '}},\quad
2^{j_2}\approx \bigl(\frac{n}{(\log
  n)^{I\{\epsilon < 0\}}}
  \bigr)^{\frac{s}{(2s+\nu')(s-\frac{1}{p}+\frac{1}{\rho})}},\quad
\lambda _j =\eta \sqrt{U_j^Y(j-j_1)/n},$$
and $\eta > 2 (\frac{2\rho \nu'}{\kappa'})^{\frac{1}{2}}$, so that we
have by Proposition $1$: $P(|\hat{\beta}_{j,k}-\beta_{j,k}|\ge
\lambda_{j}) \lesssim 2^{-\kappa' \eta^{2}(j-j_1)}.$

\bigskip

\noindent We proceed as in \cite{Donoho3} by distinguishing six terms:

\begin{align*}
\hat{f}_{n}-f&=\sum_{k \in R_{j}}
(\hat{\alpha}_{j_1,k}-\alpha_{j_1,k}) \Phi_{j,k}+ \sum _{j \ge
j_2}\sum _{k \in R_j} \beta_{j,k} \Psi_{j,k}\\&+\sum_{j_1<j<j_2}
\sum_{k \in R_j} (\hat{\beta}_{j,k}-\beta_{j,k})\Psi_{j,k}[I_{
\{|\hat{\beta}_{j,k}| \geq \lambda _j, |\beta_{j,k}|< \lambda _j
/2\}}+I_{ \{|\hat{\beta}_{j,k}| \geq \lambda _j, |\beta_{j,k}|
\geq \lambda _j /2\}}]\\&+\sum_{j_1<j<j_2} \sum_{k \in R_j}
\beta_{j,k}\Psi_{j,k}[I_{ \{|\hat{\beta}_{j,k}| < \lambda _j,
|\beta_{j,k}| \geq 2\lambda _j\}}+I_{ \{|\hat{\beta}_{j,k}| <
\lambda _j, |\beta_{j,k}| < 2 \lambda _j\}}]\\&= e_s +e_b +e_{bs}
+e_{bb} +e_{sb} +e_{ss}.\end{align*}

\noindent Like before the stochastic error is bounded by:
$$E(\|e_s\|_{\rho})\lesssim \frac{2^{\nu ' j_1}}{ \sqrt{n}}
,$$ and by using Sobolev embeddings it is easy to see
that: $$E(\|e_b\|_{\rho}) \lesssim 2^{-j_2
(s-\frac{1}{p}+\frac{1}{\rho})}.$$

The terms $e_{bs}$ and $e_{sb}$ can be grouped together because of the two following assertions: \\
$\{|\hat{\beta}_{j,k}| < \lambda _j, |\beta_{j,k}| \geq 2\lambda
_j\} \cup \{|\hat{\beta}_{j,k}| \geq \lambda _j, |\beta_{j,k}|<
\lambda _j /2\} \subset \{|\hat{\beta}_{j,k}-\beta_{j,k}|>\lambda _j
/2\}$, and \\
$[|\hat{\beta}_{j,k}| < \lambda _j, |\beta_{j,k}| \geq 2\lambda
_j] \Rightarrow  [|\beta_{j,k}|\leq
2|\hat{\beta}_{j,k}-\beta_{j,k}|]$. Consequently:

\begin{align*}
E(\|e_{bs}\|_{\rho}+\|e_{sb}\|_{\rho})&\lesssim
\sum_{j_1<j<j_2}2^{j
(\frac{1}{2}-\frac{1}{\rho})} (E\sum_{k \in R_j}
|\hat{\beta}_{j,k}-\beta_{j,k}|^{\rho}I_{\{|\hat{\beta}_{j,k}-\beta_{j,k}|>\lambda
_j /2\}})^{\frac{1}{\rho}}\\
&\le \sum_{j_1<j<j_2} 2^{j (\frac{1}{2}-\frac{1}{\rho})}
(\sum_{k \in R_j}
(E|\hat{\beta}_{j,k}-\beta_{j,k}|^{2\rho})^{\frac{1}{2}}(P\{|\hat{\beta}_{j,k}-\beta_{j,k}|>\lambda
_j /2\})^{\frac{1}{2}})^{\frac{1}{\rho}}\\
&\lesssim \sum_{j_1<j<j_2} 2^{j (\frac{1}{2}-\frac{1}{\rho})}
(\sum_{k \in R_j} \frac{2^{\rho\nu j}}{n^{\frac{\rho}{2}}}
2^{-\frac {\kappa' (\eta /2)^{2}(j-j_1)} {2}})^{\frac{1}{\rho}}\\
& \le \frac{2^{\nu' j_1}}{n^{\frac{1}{2}}}
\sum_{0<j<j_2-j_1}2^{(\nu'-\frac{\kappa'
(\eta /2)^{2}}{2\rho})j}\\
&\lesssim \frac{2^{\nu'
j_1}}{n^{\frac{1}{2}}},
\end{align*} where we used Cauchy Schwartz inequality and Proposition $1$.

\bigskip

For $e_{bb}$ we use the characterization of Besov spaces:
\begin{align*}
E(\|e_{bb}\|_{\rho})&\lesssim \sum_{j_1<j<j_2} 2^{j
(\frac{1}{2}-\frac{1}{\rho})} (\sum_{k \in R_j}
E|\hat{\beta}_{j,k}-\beta_{j,k}|^{\rho}I_{\{|\beta_{j,k}|\ge
\lambda_j /2\}}\bigr)^{\frac{1}{\rho}}\\
&\lesssim \sum_{j_1<j<j_2} 2^{j (\frac{1}{2}-\frac{1}{\rho})}
\bigl(\sum_{k \in R_j} \frac{2^{\rho \nu j}}{n^{\frac{\rho}{2}}}(\frac{|\beta_{j,k}|}{\lambda_j /2})^{p}\bigr)^{\frac{1}{\rho}}\\
&\lesssim \sum_{j_1 <j<j_2 } \bigl(\frac{2^{j (\frac{\rho}{2}-1+(\rho
-p)\nu)}}{n^{\frac{\rho-p}{2}}(j-j_1)^{\frac{p}{2}}}
2^{-pj(s+\frac{1}{2}-\frac{1}{p})}
(\|f\|^{s}_{p,\infty})^{p}\bigr)^{\frac{1}{\rho}}\\
&\lesssim \frac{1}{n^{\frac{\rho-p}{2\rho}}} \sum_{j_1 <j<j_2 } \bigl(\frac{2^{-\epsilon j}}{(j-j_1)^{\frac{p}{2}}}\bigr)^{\frac{1}{\rho}}.
\end{align*}

\bigskip

Lastly for $e_{ss}$ we remark that $|\beta_{j,k}|^{\rho}\le (2
\lambda _j)^{\rho-p}|\beta_{j,k}|^{p}$ and we use again the
characterization of Besov spaces:
\begin{align*}
E(\|e_{ss}\|_{\rho})&\lesssim \sum_{j_1<j<j_2} 2^{j
(\frac{1}{2}-\frac{1}{\rho})} \bigr( (2 \lambda _j)^{\rho-p}\sum_{k \in
R_j}|\beta_{j,k}|^{p}\bigr)^{\frac{1}{\rho}}\\
&\lesssim \sum_{j_1<j<j_2} \bigl( \frac{2^{j
(-ps+\nu'(\rho-p))}}{n^{\frac{\rho-p}{2}}}(j-j_1)^{\frac{\rho-p}{2}}(\|f\|^{s}_{p,\infty})^{p}\bigr)^{\frac{1}{\rho}}\\
&\lesssim \frac{1}{n^{\frac{\rho-p}{2\rho}}} \sum_{j_1 <j<j_2 } \bigl(2^{-\epsilon j}(j-j_1)^{\frac{\rho-p}{2}}\bigr)^{\frac{1}{\rho}}
\end{align*}

According to these bounds $e_{bs}$, $e_{sb}$ and $e_s$ are of the same order and
$e_{ss}$ dominates $e_{bb}$, so we choose $j_1$ and $j_2$ so as
to balance the bounds of $e_b$, $e_s$ and $e_{ss}$.

\noindent In the regular zone we have: $$E(\|e_{ss}\|_{\rho})\lesssim
\bigl( \frac{2^{-\epsilon j_1}}{n^{\frac{\rho-p}{2}}}\bigr)^{\frac{1}{\rho}},$$

\noindent and in the sparse zone:
$$E(\|e_{ss}\|_{\rho})\lesssim \bigl( \frac{j_2 2^{-\epsilon j_2}}{n^{\frac{\rho-p}{2}}}\bigr)^{\frac{1}{\rho}}.$$ Thus with the announced choices of $j_1$ and $j_2$ we get the prescribed
rates in both zones.

\bigskip

Lastly in the critical zone we change the majoration of
$(\beta_{j,k})$ in $e_{bb}$ and $e_{ss}$ by using:
\begin{align*}
\sum_{j_1<j<j_2} \bigr( 2^{pj(s+\frac{1}{2}-\frac{1}{p})}\sum_{k \in
R_j}|\beta_{j,k}|^{p}\bigl)^{\frac{1}{\rho}}
&\lesssim (j_2 -j_1)^{1-\frac{p}{\rho q}}(\|f\|^{s}_{p,q})^{\frac{p}{\rho}} \quad
\mbox{if} \quad \frac{p}{\rho}<q,\\
&\lesssim (\|f\|^{s}_{p,q})^{q} \quad
\mbox{if} \quad \frac{p}{\rho}\ge q.
\end{align*} Here again $e_{ss}$ is dominant and of the order: $E(\|e_{ss}\|_{\rho})\lesssim
(\frac{j_2}{n})^{\frac{\rho-p}{2\rho}}
j_2^{(1-\frac{p}{\rho q})_+},$ hence the extra logarithmic factor.

\subsubsection{Proof of the rates of the adaptive estimator}

To prove Theorem $3$ we use a theorem for thresholding algorithms
established by Kerkyacharian and Picard (Theorem $3.1$ in
\cite{maxisets}) which holds in a very general setting where one wants to
estimate an unknown function $f$ thanks to observations in a
sequence of statistical models $(E_{n})_{n \in \mathbb{N}}$. It
uses the Temlyakov inequalities, let us first recall this notion.

\begin{definition}
Let $e_{n}$ be a basis in $L^{\rho}$. It satisfies the Temlyakov
property if there are absolute constants $c$ and $C$ such that for
all $\Lambda \in \mathbb{N}$:
$$c \sum _{n \in \Lambda} \int |e_{n}(x)|^{\rho}dx \le \int \{\sum _{n \in
  \Lambda} \int |e_{n}(x)|^{2} \}^{\rho /2} dx \le C\sum _{n \in \Lambda} \int |e_{n}(x)|^{\rho}dx.$$
\end{definition}

Now let $(\psi_{j,k})_{j,k}$ denote a periodized wavelet basis and let $\rho>1$ and $0<r<\rho$. Assume that there
exist a positive value $\delta>0$, a positive sequence
$(\sigma_{j})_{j\ge -1}$, a positive sequence $c_n$ tending to
$0$, and a subset $\Lambda_{n}$ of $\mathbb{N}^{2}$ such that :

\begin{equation}
\label{h1} |\Lambda_{n}| \sim c_{n}^{-\delta} \mbox{ where }|S|
\mbox{ denotes the cardinal of the set }S,
\end{equation}

\begin{equation}
\label{h2} (\sigma_{j}\psi_{j,k})_{j,k} \mbox{ satisfies the
Temlyakov property,}
\end{equation}

\begin{equation}
\label{h3} \sup _{n} [\mu \{\Lambda_n \} c_{n}^{\rho}]< \infty,
\end{equation} where $\mu$ is the following measure on $\mathbb{N}^{2}$:
$$\mu (j,k)=\|\sigma_{j} \psi_{j,k} \|_{\rho}^{\rho}=2^{j(\rho /2-1)}\sigma_{j}^{\rho} \|\psi
\|_{\rho}^{\rho}.$$

Assume also that we have a statistical procedure yielding
estimators $\hat{\beta}_{j,k}$ of the wavelet coefficients
$\beta_{j,k}$ of $f$ in the basis $(\psi_{j,k})_{j,k}$ and a positive value $\eta>0$ such that
for all $(j,k) \in \Lambda_{n}$:

\begin{equation}
\label{h4} E(|\hat{\beta}_{j,k}-\beta_{j,k}|^{2\rho})\le C(c_n
{\sigma_j})^{2\rho},
\end{equation}

\begin{equation}
\label{h5} P(|\hat{\beta}_{j,k}-\beta_{j,k}|\ge \eta \sigma_j
c_{n}/2)\le C\min (c_{n}^{2\rho}, c_{n}^{4}).
\end{equation}

\bigskip

Finally let $l_{r,\infty}(\mu)$ and $A(c_{n}^{\rho-r})$ be the following
spaces and let $\hat{f}_{n}$ be the following estimator:
\begin{align*}
&l_{r,\infty}(\mu)=\{ f, \quad \sup _{\lambda>0}[\lambda^{q} \mu \{(j,k)/|\beta_{j,k}|>\sigma_j \lambda  \} ] < \infty\},\\
&A(c_{n}^{\rho-r})=\{ f, \quad c_{n} ^{-(\rho-r)} \| f - \sum _{\kappa \in \Lambda_{n}} \beta_{\kappa}\psi_{\kappa} \|^{\rho}_{\rho}< \infty \},\\
&\hat f _n=\sum_{j,k \in \Lambda _{n}} \hat{\beta}_{j,k} I_{ \{
  |\hat{\beta}_{j,k}| \geq \eta \sigma_{j} c_n \} } \psi_{j,k}.
\end{align*}

\begin{theorem}
Using the objects defined above and under the hypotheses
\eqref{h1} to \eqref{h5}, we have the following equivalence:
\begin{equation*}
E \|\hat{f}_{n}-f \| _{\rho}^{\rho} \lesssim c_{n}^{\rho
-r} \quad \Longleftrightarrow \quad f \in l_{r,\infty}(\mu) \cap A(c_{n}^{\rho-r}).
\end{equation*}
\end{theorem}

\bigskip

We adapt this to estimator $\hat f_n^D $ by setting, for given $\rho>1$, $p> 1$,
$s> 1/p$ and $q>1$:\\
$c_n=\sqrt{\frac{\log
    (n)^{\frac{\alpha+1}{\alpha}}}{n}}, \quad
\sigma_j=2^{\nu j}, \quad
2^{j_1}\approx\{ \frac{n} {\log
    (n)^{\frac{\alpha+1}{\alpha}}} \}^{\frac{1}{1+2\nu}}, \quad
\Lambda_{n}=\{ (j,k)\: | \: -1\le j \le j_{1}, \: k\in R_j \}.$

\bigskip

With these choices we have: $$|\Lambda_n |\asymp 2^{j_1} \asymp
c_{n}^{-2/(1+2\nu)},$$ $$\mu (\Lambda _n)=\sum _{j=0}^{j_1-1}2^{j}
2^{j(\rho /2-1)}2^{\rho \nu j}\asymp 2^{j_1 \rho(\nu+1/2)}.$$ Consequently
\eqref{h3} and \eqref{h1} hold with
$\delta=2/(1+2\nu)$. Condition \eqref{h2} is also satisfied, the
proof can be found in \cite{deconvolution}. Moreover thanks to Proposition
$1$, it is easy to establish that the estimators $\hat{\beta}_{j,k}$
used by \eqref{estimateur} satisfy \eqref{h4} and \eqref{h5} as soon
as $\eta> 2(\frac{max(2,\rho)}{\kappa})^{\frac{\alpha+1}{2\alpha}}$.

\bigskip

Then we prove Theorem $3$ by setting $r$ such that the right hand side of the
inequality in the first point of the theorem corresponds to the
rates in the sparse and in the regular case, ie:
$$r=\rho-2\rho\frac{s-1/p+1/\rho }{2s+2\nu +1-2/p},$$ or
$$r=\rho-2\rho\frac {s} {2s+2\nu +1},$$ and by showing that the space
over which the risk is maximized is included in the
maxiset, if we add the condition $q\le p$ in the critical case
$\frac{2s+2\nu +1}{\rho}=\frac{2\nu +1}{p}$:
$$M(s,p,q,R)\subset l_{r,\infty}(\mu)\cap A(c_{n}^{\rho-r}).$$ 

\bigskip

The inclusion $M(s,p,q,R)\subset A(c_{n}^{\rho-r})$ is established in
\cite{deconvolution}, and the following proof of $M(s,p,q,R)\subset
l_{r,\infty}(\mu)$ uses the same arguments as \cite{boxcar} for the
boxcar blur. We have:
\begin{align*}
\mu \{(j,k): \: |\beta_{j,k}|>2^{\nu j}\lambda \}&=\sum _{j\ge 0,\: k\in
  R_j} 2^{j (\rho(\nu+1/2)-1)}I\{ |\beta_{j,k}|>2^{\nu j}\lambda \}\\
&\le \sum _{j} (2^{j \rho(\nu+1/2)})\wedge (2^{j (\rho(\nu+1/2)-1)}\sum
  _k (|\beta_{j,k}|/(2^{\nu j}\lambda))^{p}\\
&\le \sum _{j} (2^{j \rho(\nu+1/2)})\wedge
  (\frac{2^{-j(sp+\nu'p-\nu'\rho)}}{\lambda^p}\epsilon_j ^{p}),
\end{align*} where $\nu'=\nu+1/2$ and $\epsilon_j \in l_q$. We cut the
  sum at $J$ such that $2^J \asymp \lambda^{-r/(\nu' \rho)}.$

\bigskip

In the regular case we have:
\begin{align*}
\mu \{(j,k): \: |\beta_{j,k}|>2^{\nu j}\lambda \}&\le \lambda
^{-r}+\frac{\lambda^{(sp-\nu'(\rho-p))\frac{r}{\nu' \rho}}}{\lambda^p},
\end{align*} and the power of $\lambda$ in the second term is also
exactly $-r$.

\bigskip

In the critical case we obtain, since $q \le p$:
\begin{align*}
\mu \{(j,k): \: |\beta_{j,k}|>2^{\nu j}\lambda \}&\le \lambda
^{-r}+\frac{\sum_j \epsilon_j^p}{\lambda^p}\lesssim \lambda^{-r}+\frac{\sum_j \epsilon_j^q}{\lambda^p}\lesssim \lambda^{-r}+\lambda^{-p},
\end{align*} and $r=p$ in this case.

\bigskip

Lastly in the sparse case (where $r\ge p$ is satisfied) we
use the Sobolev embedding $B^s_{p,q} \subset B^{s'}_{r,q}$
with $s'=s-1/p+1/r$. We proceed as before by cutting the sum at $J$
such that $2^J \asymp \lambda^{-r/(\nu' \rho)} $ and noticing that
$s'r+\nu'r-\nu'\rho=0$. There exists $\tilde{\epsilon}_j \in l_r$ such
that:
\begin{align*}
\mu \{(j,k): \: |\beta_{j,k}|>2^{\nu j}\lambda \}
&\le \sum _{j}(2^{j \rho\nu'})\wedge (2^{j (\rho\nu'-1)}\sum
  _k (|\beta_{j,k}|/(2^{\nu j}\lambda))^{r}\\
&\le \sum _{j} (2^{j \rho\nu'})\wedge
  (\frac{\tilde{\epsilon}_j ^{r}}{\lambda^r})\\
&\lesssim \lambda^{-r}.
\end{align*}

Thus $\mu \{(j,k): \: 2^{\nu j}\lambda \}\lesssim 1/\lambda^{r}$ for
both values of $r$, and finally using the equivalence in Theorem $4$ and Jensen inequality we obtain the
prescribed rates for $E \|\hat{f}_{n}^{D}-f \| _{\rho}$.

\bigskip
\noindent
\textbf{Acknowledgments}

\bigskip
\noindent
I would like to thank Professor Dominique Picard for introducing me to
the subject and for her numerous and helpful comments.

\appendix
\newpage

\bibliographystyle{abbrvnat}
\bibliography{dec_random_filter}

\begin{thebibliography}{33}
\expandafter\ifx\csname natexlab\endcsname\relax\def\natexlab#1{#1}\fi
\expandafter\ifx\csname url\endcsname\relax
  \def\url#1{{\tt #1}}\fi

\bibitem[Abramovich and Silverman(1998)]{Silverman}
F.~Abramovich and B.~Silverman.
\newblock Wavelet decomposition approaches to statistical inverse problems.
\newblock {\em Biometrika}, 85\penalty0 (1):\penalty0 115--129, 1998.

\bibitem[Bertero and Boccacci(1998)]{Bertero}
M.~Bertero and P.~Boccacci.
\newblock Introduction to inverse problems in imaging.
\newblock {\em Institute of Physics, Bristol and Philadelphia}, 1998.

\bibitem[Butucea(2004)]{butucea}
C.~Butucea.
\newblock Deconvolution of supersmooth densities with smooth noise.
\newblock {\em Canad. J. Statist.}, 32\penalty0 (2):\penalty0 181--192, 2004.

\bibitem[Cavalier et~al.(2003)Cavalier, Golubev, Lepski, and
  Tsybakov]{tsybakov1}
L.~Cavalier, Y.~Golubev, O.~Lepski, and A.~Tsybakov.
\newblock Block thresholding and sharp adaptive estimation in severely
  ill-posed inverse problems.
\newblock {\em Theory of Probability and its Applications}, 48\penalty0
  (3):\penalty0 534--556, 2003.

\bibitem[Cavalier and Hengartner(2004)]{cavalier}
L.~Cavalier and N.~W. Hengartner.
\newblock Adaptive estimation for inverse problems with noisy operators.
\newblock {\em Manuscript}, 2004.

\bibitem[Cavalier and Tsybakov(2002)]{cavalier_tsybakov}
L.~Cavalier and A.~Tsybakov.
\newblock Sharp adaptation for inverse problems with random noise.
\newblock {\em Probab. Theory Related Fields}, 123\penalty0 (3):\penalty0
  323--354, 2002.

\bibitem[Cohen et~al.(2002)Cohen, Hoffmann, and Reiss]{marc}
A.~Cohen, M.~Hoffmann, and M.~Reiss.
\newblock Adaptive wavelet galerkin methods for linear inverse problems.
\newblock {\em Preprint LPMA}, dec 2002.

\bibitem[Cohen et~al.(2004)Cohen, Hoffmann, and Reiss]{hoffmann}
A.~Cohen, M.~Hoffmann, and M.~Reiss.
\newblock On adaptive estimation in linear inverse problems with error in the
  operator.
\newblock {\em Manuscript}, 2004.

\bibitem[Donoho(1995)]{Donoho95}
D.~Donoho.
\newblock Nonlinear solution of linear inverse problems by wavelet-vaguelette
  decomposition.
\newblock {\em Applied Computational and Harmonic Analysis}, 2:\penalty0
  101--126, 1995.

\bibitem[Donoho and Johnstone(1994)]{dj94}
D.~Donoho and I.~Johnstone.
\newblock Ideal spatial adaptation by wavelet shrinkage.
\newblock {\em Biometrika}, 81:\penalty0 425--455, 1994.

\bibitem[Donoho et~al.(1996)Donoho, Johnstone, Kerkyacharian, and
  Picard]{Donoho3}
D.~Donoho, I.~Johnstone, G.~Kerkyacharian, and D.~Picard.
\newblock Density estimation by wavelet thresholding.
\newblock {\em Ann. Statist.}, 2:\penalty0 508--539, 1996.

\bibitem[Donoho et~al.(1997)Donoho, Johnstone, Kerkyacharian, and
  Picard]{Donoho2}
D.~Donoho, I.~Johnstone, G.~Kerkyacharian, and D.~Picard.
\newblock Universal near minimaxity of wavelet shrinkage.
\newblock {\em Festschrift for Lucien Le Cam}, pages 183--218, 1997.

\bibitem[Efromovich and Koltchinskii(2001)]{efro}
S.~Efromovich and V.~Koltchinskii.
\newblock On inverse problems with unknown operators.
\newblock {\em IEEE}, 47\penalty0 (7):\penalty0 2876--2894, 2001.

\bibitem[Fan and Koo(2002)]{Fan}
J.~Fan and J.~Koo.
\newblock Wavelet deconvolution.
\newblock {\em IEEE Transactions on Information Theory}, 48\penalty0
  (3):\penalty0 734--747, 2002.

\bibitem[H\"ardle et~al.(1998)H\"ardle, Kerkyacharian, Picard, and
  Tsybakov]{wavelets}
W.~H\"ardle, G.~Kerkyacharian, D.~Picard, and A.~Tsybakov.
\newblock {\em Wavelets, Approximation and Statistical Applications}.
\newblock Springer-Verlag, 1998.

\bibitem[Harsdorf and Reuter(2000)]{Harsdorf}
S.~Harsdorf and R.~Reuter.
\newblock Stable deconvolution of noisy lidar signals.
\newblock {\em Submitted to Earsel meeting}, 2000.

\bibitem[Johnstone et~al.(2004)Johnstone, Kerkyacharian, Picard, and
  Raimondo]{deconvolution}
I.~M. Johnstone, G.~Kerkyacharian, D.~Picard, and M.~Raimondo.
\newblock Wavelet deconvolution in a periodic setting.
\newblock {\em Journal of the Royal Statistical Society}, 66\penalty0
  (3):\penalty0 1--27, 2004.

\bibitem[Johnstone and Raimondo(2004)]{Johnstone}
I.~M. Johnstone and M.~Raimondo.
\newblock Periodic boxcar deconvolution and diophantine approximation.
\newblock {\em Annals of Statistics}, 32\penalty0 (5):\penalty0 1781--1804,
  2004.

\bibitem[Kalifa and Mallat(2003)]{Kalifa}
J.~Kalifa and S.~Mallat.
\newblock Thresholding estimators for linear inverse problems and
  deconvolutions.
\newblock {\em Annals of Statistics}, 31:\penalty0 58--109, 2003.

\bibitem[Kerkyacharian and Picard(2000)]{maxisets}
G.~Kerkyacharian and D.~Picard.
\newblock Thresholding algorithms, maxisets and well concentrated bases.
\newblock {\em Test}, 9\penalty0 (2), 2000.

\bibitem[Kerkyacharian et~al.(2004)Kerkyacharian, Picard, and Raimondo]{boxcar}
G.~Kerkyacharian, D.~Picard, and M.~Raimondo.
\newblock Adaptive boxcar deconvolution on full lebesgue measure sets.
\newblock {\em Preprint LPMA}, 2004.

\bibitem[Korostelev and Tsybakov(1993)]{korostelev_tsybakov93}
V.~Korostelev and A.~Tsybakov.
\newblock {\em Minimax theory of image reconstruction}.
\newblock Springer-Verlag, 1993.

\bibitem[Mallat(1998)]{Mallat}
S.~Mallat.
\newblock {\em A wavelet tour of signal processing (2nd edition)}.
\newblock Academic Press Inc., San Diego, CA, 1998.

\bibitem[Meyer(1990)]{Meyer}
Y.~Meyer.
\newblock {\em Ondelettes et Op\'erateurs-I}.
\newblock Hermann, 1990.

\bibitem[Nussbaum and Pereverzev(1999)]{illpose}
M.~Nussbaum and S.~Pereverzev.
\newblock The degree of ill-posedness in stochastic and deterministic noise
  models.
\newblock {\em Preprint WIAS}, 1999.

\bibitem[OFTA(1999)]{ofta}
OFTA.
\newblock {\em Probl\`emes inverses : de l'exp\'erimentation \`a la
  mod\'elisation}.
\newblock Observatoire Francais des Techniques Avanc\'ees, 1999.

\bibitem[Pensky and Vidakovic(1999)]{Pensky}
M.~Pensky and B.~Vidakovic.
\newblock Adaptive wavelet estimator for nonparametric density deconvolution.
\newblock {\em Annals of Statistics}, 27:\penalty0 2033--2053, 1999.

\bibitem[Reiss(2001)]{markus2}
M.~Reiss.
\newblock Nonparametric estimation for stochastic delay differential equations.
\newblock {\em Ph.D. thesis, Humboldt Universit\"at zu Berlin}, 2001.

\bibitem[Reiss(2004)]{markus}
M.~Reiss.
\newblock Adaptive estimation for affine stochastic delay differential
  equations.
\newblock {\em Submitted to Bernoulli}, 2004.

\bibitem[Tsybakov(2000)]{tsybakov3}
A.~Tsybakov.
\newblock On the best rate of adaptive estimation in some inverse problems.
\newblock {\em C.R. Acad. Sci.}, 1\penalty0 (330):\penalty0 835--840, 2000.

\bibitem[Tsybakov(2001)]{tsybakov4}
A.~Tsybakov.
\newblock Sharp adaptive estimation of linear functionals.
\newblock {\em Ann. Statist.}, 29\penalty0 (6):\penalty0 1567--1600, 2001.

\bibitem[Tsybakov(2004)]{tsybakov03}
A.~Tsybakov.
\newblock {\em Introduction \`a l'estimation Non-param\'etrique}.
\newblock Springer, 2004.

\bibitem[Walter and Shen(1999)]{Walter}
G.~Walter and X.~Shen.
\newblock Deconvolution using the meyer wavelet.
\newblock {\em Journal of Integral Equations and Applications}, 11:\penalty0
  515--534, 1999.

\end{thebibliography}

\end{document}